\documentclass[12pt]{amsart}
\usepackage{graphicx}
\usepackage{amssymb}
\usepackage{amsfonts}
\usepackage{amsmath,amscd}
\usepackage{array}
\usepackage{rotating}
\usepackage{epsfig}

\newtheorem{example}{Example}[section]

\newtheorem{theorem}[example]{Theorem}

\newtheorem{corollary}[example]{Corollary}

\newtheorem{proposition}[example]{Proposition}

\def\Proof{\noindent \it Proof -- \rm}
\def\qed{\hspace{3.5mm} \hfill \vbox{\hrule height 3pt depth 2 pt width 2mm}
\bigskip}

\def\QSym{{\it QSym}}          % QSym
\def\NCSF{{\bf Sym}}           % NCSF
\def\FQSym{{\bf FQSym}}        % permutations
\def\FSym{{\bf FSym}}          % tableaux standard
\def\PBT{{\bf PBT}}            % arbres binaires planaires
\def\PQSym{{\bf PQSym}}        % Parking 
\def\WQSym{{\bf WQSym}}        % Mots initiaux 
\def\CQSym{{\bf CQSym}}        % Catalan
\def\CPQSym{{\it PQSym}}       % Parking (commutative)
\def\EQSym{{\it EQSym}}        % Endofonctions
\def\ESym{{\bf ESym}}          % Endofonctions duales
\def\PiQSym{{\it {\Pi}QSym}}   % Mots initiaux
\def\CCQSym{{\it CQSym}}       % Parking croissantes (commutative)
\def\SGQSym{{\it \SG QSym}}    % Cycles
\def\SGSym{{\bf \SG Sym}}      % Cycles duale
\def\PhiSym{{\bf {\Phi}Sym}}   % Cycles bis (duale precedente)

\def\WSym{{\bf WSym}}          % Mots symetriques

\def\MEF{{\it M}}              % base naturelle des endofonctions (EQSym)
\def\SEF{{\bf S}}              % base naturelle des endofonctions duales
\def\uq{{\it U}}               % base naturelle des compositions (QSym)
\def\upi{{\it U}}              % base naturelle des parts ensembl. (PiQSym)
\def\ul{{\it u}}               % base des images des U_\pi (fns syms)
\def\Mper{{\it M}}             % base naturelle des cycles (SGQSym)
\def\Sper{{\bf S}}             % base naturelle des cycles duaux
\def\Mpa{{\it M}}              % base naturelle des parking commut. (PQSym)
\def\Mw{{\bf M}}               % base de WSym (mots par action symetrique)

\def\csupp{{\rm supp}}         % support = $\pi$ des cycles de perm.
\def\picyc{{\shuffle}}         % produit "shuffle" de deux cycles (PhiSym)
\def\picycgen{{\smile}}        % produit "Wick de 2 cycles (PhiSym)

\def\ncbinomial#1#2{\left[\,\begin{matrix}#1 \cr #2\end{matrix}\,\right]}
\def\myphi{{\bf {\bf\phi}}} % base des cycles bis

\def\inv{{\rm inv}}     % inversions

\def\sconc{\bullet}     % concatenation shiftee
\def\Std{{\rm Std}}     % standardisation
\def\cstd{{\rm cstd}}   % standardisation cyclique
\def\Park{{\rm Park}}   % parkisation

\def\<{\langle}
\def\>{\rangle}
\def\ZZ{{\mathbb Z}}    % entiers relatifs
\def\CC{{\mathbb C}}    % complexes
\def\KK{{\mathbb K}\, } % corps K

\def\park{{\bf a}} % les parking 
\def\parkc{{\pi}} % les parking croissantes
\def\tr{\operatorname{tr}}

\def\F{{\bf F}}         % F de FQSym
\def\S{{\bf S}}         % S de NCSF
\def\G{{\bf G}}         % G de FQSym^*
\def\SG{{\mathfrak S}}  % groupe symetrique
\def\gr{{\rm gr}}

\def\goth{\mathfrak}
\def\tensor{\otimes}
\def\dim{{\rm dim}}

\def\ch{\operatorname {ch}}

\def\End{\operatorname{End}} % endomorphismes

\def\PF{{\rm PF}}   % fonctions de parking
\def\shuff#1#2{\mathbin{
\hbox{\vbox{ \hbox{\vrule \hskip#2 \vrule height#1 width 0pt
}%
\hrule}%
\vbox{ \hbox{\vrule \hskip#2 \vrule height#1 width 0pt
\vrule }%
\hrule}%
}}}

\long\def\psboxit#1#2{%
\begingroup\setbox0=\hbox{#2}%
\dimen0=\ht0 \advance\dimen0 by \dp0%
    \hbox{%
    \copy0%
    }%hbox
\endgroup%
}%psboxit
\def\SetTableau#1#2#3#4{%
  \gdef\Tabvrule{\vrule\vrule width-0.4pt}
  \gdef\Tabhrule{\hrule\hrule height-0.4pt}  
  \gdef\Tabstrut{\vrule height#1 depth#2 width0pt\relax}
  \gdef\Tabbox##1{\hbox to #3{\hskip0.4pt\hfill\Tabstrut$#4##1$\hfill}}
} %setTableau
\def\TasseTableau{\SetTableau{1.48ex}{0.32ex}{1.8ex}{\scriptstyle}}
\def\PetitTableau{\SetTableau{1.65ex}{0.55ex}{2.2ex}{\scriptstyle}}

\def\Case#1{\vcenter{\Tabhrule%
                   \hbox{\Tabvrule\Tabbox{#1}\Tabvrule}\Tabhrule}}
\def\GenTab#1{\vcenter{\halign{&$\Case{##}$\cr#1}}\egroup}
\def\Tableau{%         
  \bgroup%
  \let\ =\omit%
  \let\\=\cr%
  \offinterlineskip\GenTab}

\def\qbin#1#2{\begin{bmatrix} #1 \\ #2\end{bmatrix}}

\def\shuf{{\mathchoice{\shuff{7pt}{3.5pt}}%
{\shuff{6pt}{3pt}}%
{\shuff{4pt}{2pt}}%
{\shuff{3pt}{1.5pt}}}}%
\def\shuffle{\,\shuf\,}

\def\SS{\mathbf S} % Base FSym
\def\Pp{{\bf P}}        % P de PBT

\title[Commutative combinatorial Hopf algebras]%
{Commutative combinatorial Hopf algebras}
\author[F. Hivert, J.-C.~Novelli, and J.-Y.~Thibon]%
{Florent Hivert, Jean-Christophe Novelli, and Jean-Yves Thibon}

\address[] {Institut Gaspard Monge, Universit\'e de Marne-la-Vall\'ee \\
5 Boulevard Descartes \\Champs-sur-Marne \\77454 Marne-la-Vall\'ee cedex 2 \\
FRANCE}
\email[Florent Hivert]{hivert@univ-mlv.fr}
\email[Jean-Christophe Novelli]{novelli@univ-mlv.fr}
\email[Jean-Yves Thibon]{jyt@univ-mlv.fr} 
\date{}

\begin{document}

\begin{abstract}
We propose several constructions of commutative or cocommutative Hopf algebras
based on various combinatorial structures, and investigate the relations
between them. A commutative Hopf algebra of permutations is obtained by a
general construction based on graphs, and its non-commutative dual is realized
in three different ways, in particular as the Grossman-Larson algebra of heap
ordered trees.
Extensions to endofunctions, parking functions, set compositions, set
partitions, planar binary trees and rooted forests are discussed. Finally, we
introduce one-parameter families interpolating between different structures
constructed on the same combinatorial objects.
\end{abstract}

\maketitle

{\footnotesize
\tableofcontents
}

%%%%%%%%%%%%%%%%%%%%%%%%%%%%%%%%%%%%%%%%%%%%%%%%%%%%%%%%%%%%%%%%%%%%%%%%%%%%%%%
%%%%%%%%%%%%%%%%%%%%%%%%%%%%%%%%%%%%%%%%%%%%%%%%%%%%%%%%%%%%%%%%%%%%%%%%%%%%%%%
%%%%%%%%%%%%%%%%%%%%%%%%%%%%%%%%%%%%%%%%%%%%%%%%%%%%%%%%%%%%%%%%%%%%%%%%%%%%%%%
\section{Introduction}

Many examples of Hopf algebras based on combinatorial structures are known.
Among these, algebras based on permutations and planar binary trees play a
prominent role, and arise in seemingly unrelated
contexts~\cite{MR,LR1,NCSF6,BF}.
As Hopf algebras, both are noncommutative and non cocommutative, and in fact
self-dual.

More recently, cocommutative Hopf algebras of binary trees and permutations
have been constructed~\cite{NT1,AS}. In~\cite{NT1}, binary trees arise as sums
over rearrangements classes in an algebra of parking functions, while
in~\cite{AS}, cocommutative Hopf algebras are obtained as the graded
coalgebras associated with coradical fitrations.

In~\cite{NTT}, a general method for constructing commutative Hopf algebras
based on various kind of graphs has been presented.  The aim of this note is
to investigate Hopf algebras based on permutations and trees constructed by
the method developed in~\cite{NTT}. These commutative algebras are, by
definition, realized in terms of polynomials in an infinite set of doubly
indexed indeterminates. The dual Hopf algebras are then realized by means of
non commutative polynomials in variables $a_{ij}$. We show that these first
resulting algebras are isomorphic (in a non trivial way) to the duals of those
of~\cite{AS}, and study some generalizations such as endofunctions, parking
functions, set partitions, trees, forests, and so on.

The possibility to obtain in an almost systematic way commutative, and in
general non cocommutative, versions of the usual combinatorial Hopf algebras
leads us to conjecture that these standard versions should be considered as
some kind of quantum groups, \emph{i.e.}, can be incorporated into
one-parameter families, containing an enveloping algebra and its dual for
special values of the parameter.
A few results supporting this point of view are presented in the final
section.

\medskip
In all the paper, $\KK$ will denote a field of characteristic zero.
All the notations used here is as in~\cite{NCSF1, NTT}.
This paper is an expanded and updated version of the preprint~\cite{HNT-co}.

%%%%%%%%%%%%%%%%%%%%%%%%%%%%%%%%%%%%%%%%%%%%%%%%%%%%%%%%%%%%%%%%%%%%%%%%%%%%%%%
%%%%%%%%%%%%%%%%%%%%%%%%%%%%%%%%%%%%%%%%%%%%%%%%%%%%%%%%%%%%%%%%%%%%%%%%%%%%%%%
%%%%%%%%%%%%%%%%%%%%%%%%%%%%%%%%%%%%%%%%%%%%%%%%%%%%%%%%%%%%%%%%%%%%%%%%%%%%%%%
\section{A commutative Hopf algebra of endofunctions}
\label{eqsym}

Permutations can be regarded in an obvious way as labelled and oriented graphs
whose connected components are cycles.
Actually, arbitrary \emph{endofunctions} (functions from $[n]:=\{1,\ldots,n\}$
to itself) can be regarded as labelled graphs, connecting $i$ with $f(i)$ for
all $i$ so as to fit in the framework of~\cite{NTT}, where a general process
for building Hopf algebras of graphs is described.

In the sequel, we identify an endofunction $f$ of $[n]$ with the word
\begin{equation}
w_f=f(1)f(2)\cdots f(n) \in [n]^n.
\end{equation}

Let $\{x_{i\,j} \,|\, i,j\geq1\}$ be an infinite set of commuting
indeterminates, and let $\mathcal J$ be the ideal of
$R=\KK[x_{i\,j} \,|\, i,j\geq1]$ generated by the relations
\begin{equation}
x_{i\,j} x_{i\,k}=0
\quad \text{for all $i,j,k$.}
\end{equation}
For an endofunction $f:[n]\to [n]$, define
\begin{equation}
\label{mf-eqs}
\MEF_f := \sum_{i_1 < \cdots < i_n}
            x_{i_1\, i_{f(1)}}\cdots x_{i_n\, i_{f(n)}},
\end{equation}
in $R/{\mathcal J}$.

It follows from~\cite{NTT}, Section~4, that

\begin{theorem}
The $\MEF_f$ span a subalgebra $\EQSym$ of the commutative algebra
$R/{\mathcal J}$.
More precisely, there exist non-negative integers $C_{f,g}^{h}$
such that
\begin{equation}
\MEF_f \MEF_g = \sum_{h} C_{f,g}^{h} \MEF_h.
\end{equation}
\qed
\end{theorem}

\begin{example}
{\rm
\begin{equation}
\MEF_{1} \MEF_{22} = \MEF_{133} + \MEF_{323} + \MEF_{223}.
\end{equation}
\begin{equation}
\MEF_{1} \MEF_{331} = \MEF_{1442} + \MEF_{4241} + \MEF_{4431} + \MEF_{3314}.
\end{equation}
\begin{equation}
\MEF_{12} \MEF_{21} = \MEF_{1243} + \MEF_{1432} + \MEF_{4231} + \MEF_{1324} +
\MEF_{3214} + \MEF_{2134}.
\end{equation}
\begin{equation}
\label{mf12-22}
\MEF_{12} \MEF_{22} = \MEF_{1244} + \MEF_{1434} + \MEF_{4234} + \MEF_{1334} +
\MEF_{3234} + \MEF_{2234}.
\end{equation}
\begin{equation}
\MEF_{12} \MEF_{133} = 3\MEF_{12355} + 2\MEF_{12445} + 2\MEF_{12545}
+ \MEF_{13345} + \MEF_{14345} + \MEF_{15345}.
\end{equation}
}
\end{example}

The \emph{shifted concatenation} of two endofunctions $f:[n]\to[n]$ and
$g:[m]\to[m]$ is the endofunction
$h:=f\bullet g$ of $[n+m]$ such that $w_h :=w_f\bullet w_g$, that is
\begin{equation}
\left\{
\begin{aligned}
h(i)=f(i)     \ & \text{if} &\ i\leq n\\ 
h(i)=n+g(i-n) \ & \text{if} &\ i > n\\ 
\end{aligned}
\right.
\end{equation}

We can now give a combinatorial interpretation of the coefficient
$C_{f,g}^{h}$: if $f:[n]\to[n]$ and $g:[m]\to[m]$, this coefficient is the
number of permutations $\tau$ in the shuffle product
$(1\ldots p)\shuffle (p+1\ldots p+n)$ such that
\begin{equation}
h = 
\tau^{-1} \circ (f\bullet g) \circ \tau.
\end{equation}
For example, with $f=12$ and $g=22$, one finds the set (see
Equation~(\ref{mf12-22}))
\begin{equation}
\label{mf12-22a}
\{1244, 1434, 4234, 1334, 3234, 2234 \}.
\end{equation}

Now, still following~\cite{NTT}, define a coproduct by

\begin{equation}
\label{coprodM-endof}
\Delta \MEF_h := \sum_{(f,g) ; f\sconc g=h} 
\MEF_{f} \tensor \MEF_{g}.
\end{equation}
This endows $\EQSym$ with a (commutative, non cocommutative) Hopf algebra
structure.

\begin{example}
\begin{equation}
\Delta \MEF_{626124}  = \MEF_{626124} \tensor 1 + 1\tensor \MEF_{626124}.
\end{equation}
\begin{equation}
\Delta \MEF_{4232277} = \MEF_{4232277} \tensor 1 +
                        \MEF_{42322}\tensor \MEF_{22} +
                        1\tensor \MEF_{4232277}.
\end{equation}
\end{example}

Define a \emph{connected} endofunction as a function that cannot be obtained
by non trivial shifted concatenation.
For example, the connected endofunctions for $n=1,\ 2,\ 3$ are
\begin{equation}
\begin{split}
& 1, \qquad\qquad 11,\ 21,\ 22, \\
& 111,\ 112,\ 121,\ 131,\ 211,\ 212,\ 221,\ 222,\ 231,\ 232,\\
& 233,\ 311,\ 312,\ 313,\ 321,\ 322,\ 323,\ 331,\ 332,\ 333,\\
\end{split}
\end{equation}
and the generating series of their number begins with
\begin{equation}
 t + 3\,t^2 + 20\,t^3 + 197\,t^4 + 2511\,t^5 + 38924\,t^6 +
 708105\,t^7 + 14769175\,t^8\,.
\end{equation}

Then, the definition of the coproduct of the $\MEF_f$ implies

\begin{proposition}
If $(\SEF^f)$ denotes the dual basis of $(\MEF_f)$, the graded
dual $\ESym := \EQSym^*$ is free over the set
\begin{equation}
\{ \SEF^f \,|\, f \text{\ connected} \}.
\end{equation}
\end{proposition}

Indeed, Equation~(\ref{coprodM-endof}) is equivalent to
\begin{equation}
\SEF^f \SEF^g = \SEF^{f\sconc g}.
\end{equation}
\qed

Now, $\ESym$ being a graded connected cocommutative Hopf algebra, it follows
from the Cartier-Milnor-Moore theorem that
\begin{equation}
\ESym = U(L),
\end{equation}
where $L$ is the Lie algebra of its primitive elements.
Let us now prove
\begin{theorem}
\label{prim-lib}
As a graded Lie algebra, the primitive Lie algebra $L$ of $\ESym$ is free over
a set indexed by connected endofunctions.
\end{theorem}

\Proof
Assume it is the case. By standard arguments on generating series, one finds
that the number of generators of $L$ in degree $n$ is equal to the number of
algebraic generators of $\ESym$ in degree $n$, parametrized for example by
connected endofunctions.
% (series beginning by
%$(1,3,20,197,2511,38924,\ldots)$).
We will now show that $L$ has at least this number of generators and that
those generators are algebraically independent, determining completely the
dimensions of the homogeneous components $L_n$ of $L$ whose generating series
begins by
\begin{equation}
t+3\,t^2 + 23\,t^3 + 223\,t^4 + 2800\,t^5 + 42576\,t^6 +
763220\,t^7 + 15734388\,t^8 + \ldots
\end{equation}
Following Reutenauer~\cite{Reu} p.~58, denote by $\pi_1$ the Eulerian
idempotent, that is, the endomorphism of $\ESym$ defined by
$\pi_1=\log^*(Id)$. It is obvious, thanks to the definition of $\S^f$ that
\begin{equation}
\pi_1(\S^{f}) = \S^f + \cdots,
\end{equation}
where the dots stand for terms $\SEF^g$ such that $g$ is not connected.
Since the $\SEF^f$ associated with connected endofunctions are
independent, the dimension of $L_n$ is at least equal to the
number of connected endofunctions of size $n$. So $L$ is free over a set of
primitive elements parametrized by connected endofunctions.
\qed

%%%%%%%%%%%%%%%%%%%%%%%%%%%%%%%%%%%%%%%%%%%%%%%%%%%%%%%%%%%%%%%%%%%%%%%%%%%%%%%

There are many Hopf subalgebras of $\EQSym$ which can be defined by imposing
natural restrictions to maps: being bijective (see Section~\ref{sgqsym}),
idempotent ($f^2=f$), involutive ($f^2=id$), or more generally the Burnside
classes ($f^p=f^q$), and so on.
We shall start with the Hopf algebra of permutations.

%%%%%%%%%%%%%%%%%%%%%%%%%%%%%%%%%%%%%%%%%%%%%%%%%%%%%%%%%%%%%%%%%%%%%%%%%%%%%%%
%%%%%%%%%%%%%%%%%%%%%%%%%%%%%%%%%%%%%%%%%%%%%%%%%%%%%%%%%%%%%%%%%%%%%%%%%%%%%%%
%%%%%%%%%%%%%%%%%%%%%%%%%%%%%%%%%%%%%%%%%%%%%%%%%%%%%%%%%%%%%%%%%%%%%%%%%%%%%%%
\section{A commutative Hopf algebra of permutations}
\label{sgqsym}

\medskip
We will use two different notations for permutations depending whether they
are considered as bijections from $[1,n]$ onto itself or as product of cycles.
In the first case, we will write $\sigma=31542$ for the bijection $\sigma$
where $\sigma(i)=\sigma_i$.
In the second case, the same permutation will be written $\sigma=(1352)(4)$
since $31542$ is composed of two cycles : the cycle $(1352)$ sending each
element to the next one (circularly) in the sequence and the cycle $(4)$
composed of only one element.

%%%%%%%%%%%%%%%%%%%%%%%%%%%%%%%%%%%%%%%%%%%%%%%%%%%%%%%%%%%%%%%%%%%%%%%%%%%%%%%
\subsection{The Hopf algebra of bijective endofunctions}

Let us define $\SGQSym$ as the subalgebra of $\EQSym$ spanned by the
\begin{equation}
\Mper_\sigma = \sum_{i_1 < \cdots < i_n}
            x_{i_1\, i_{\sigma(1)}}\cdots x_{i_n\, i_{\sigma(n)}},
\end{equation}
where $\sigma$ runs over bijective endofunctions, \emph{i.e.}, permutations.
Note that $\SGQSym$ is also isomorphic to the image of $\EQSym$ in the
quotient of $R/{\mathcal J}$ by the relations

\begin{equation}
x_{i\,k}x_{j\,k}=0 \quad \text{for all $i,j,k$.}
\end{equation}

By the usual argument, it follows that

\begin{proposition}
The $\Mper_\sigma$ span a Hopf subalgebra $\SGQSym$ of the commutative Hopf
algebra $\EQSym$.
\qed
\end{proposition}

As already mentionned, there exist non-negative integers
$C_{\alpha,\beta}^{\gamma}$ such that
\begin{equation}
\label{defC}
\Mper_\alpha \Mper_\beta = \sum_{\gamma} C_{\alpha,\beta}^{\gamma}
\Mper_\gamma.
\end{equation}

The combinatorial interpretation of the coefficients $C_{f,g}^h$ seen in
Section~\ref{eqsym} can be reformulated in the special case of permutations.
Write $\alpha$ and $\beta$ as a union of disjoint cycles. Split the set
$[n+m]$ into a set $A$ of $n$ elements, and its complement $B$, in
all possible ways. For each splitting, apply to $\alpha$ (resp. $\beta$) in
$A$ (resp. $B$) the unique increasing morphism of alphabets from
$[n]$ to $A$ (resp. from $[m]$ to $B$) and return the list of permutations
with the resulting cycles.
On the example $\alpha=(1)(2)=12$ and $\beta=(13)(2)=321$, this yields
\begin{equation}
\label{12-321b}
\begin{split}
(1)(2)(53)(4),\ (1)(3)(52)(4),\ (1)(4)(52)(3),\ (1)(5)(42)(3),\ (2)(3)(51)(4),\\
(2)(4)(51)(3),\ (2)(5)(41)(3),\ (3)(4)(51)(2),\ (3)(5)(41)(2),\ (4)(5)(31)(2).
\end{split}
\end{equation}
This set corresponds to the permutations and multiplicities of
Equation~(\ref{12-321c}).

A third interpretation of this product comes from the dual coproduct point of
view: $C_{\alpha,\beta}^{\gamma}$ is the number of ways of getting
$(\alpha,\beta)$ as the standardized words of pairs $(a,b)$ of two
complementary subsets of cycles of $\gamma$.
For example, with $\alpha=12$, $\beta=321$, and $\gamma=52341$, one has three
solutions for the pair $(a,b)$, namely
\begin{equation}
((2)(3), (4)(51)),\ \  ((2)(4), (3)(51)),\ \ ((3)(4), (2)(51)),
\end{equation}
which is coherent with Equations~(\ref{12-321b}) and (\ref{12-321c}).

\begin{example}
\begin{equation}
\Mper_{12\cdots n} \Mper_{12\cdots p} =
  \binom{n+p}{n} \Mper_{12\cdots (n+p)}.
\end{equation}
\begin{equation}
\Mper_{1} \Mper_{21} = \Mper_{132} + \Mper_{213} + \Mper_{321}.
\end{equation}
\begin{equation}
\Mper_{12} \Mper_{21} = \Mper_{1243} + \Mper_{1324} + \Mper_{1432} +
\Mper_{2134} + \Mper_{3214} + \Mper_{4231}.
\end{equation}
\begin{equation}
\label{12-321c}
\Mper_{12} \Mper_{321} = \Mper_{12543} + \Mper_{14325} + 2 \Mper_{15342} +
\Mper_{32145} + 2 \Mper_{42315} + 3 \Mper_{52341}.
\end{equation}
\begin{equation}
\begin{split}
\Mper_{21} \Mper_{123} &=\
      \Mper_{12354} + \Mper_{12435} + \Mper_{12543} + \Mper_{13245} +
      \Mper_{14325}\\
  &+\ \Mper_{15342} + \Mper_{21345} + \Mper_{32145} + \Mper_{42315} +
      \Mper_{52341}.
\end{split}
\end{equation}
\begin{equation}
\begin{split}
\Mper_{21} \Mper_{231} &=\
   \Mper_{21453} + \Mper_{23154} + \Mper_{24513} + \Mper_{25431} +
   \Mper_{34152}\\
 &+\ \Mper_{34521} + \Mper_{35412} + \Mper_{43251} + 
   \Mper_{43512} + \Mper_{53421}.
\end{split}
\end{equation}

\end{example}

%%%%%%%%%%%%%%%%%%%%%%%%%%%%%%%%%%%%%%%%%%%%%%%%%%%%%%%%%%%%%%%%%%%%%%%%%%%%%%%
\subsection{Duality}

Recall that the coproduct is given by

\begin{equation}
\label{coprodM}
\Delta \Mper_\sigma := \sum_{(\alpha,\beta) ; \alpha\sconc \beta=\sigma} 
\Mper_{\alpha} \tensor \Mper_{\beta}.
\end{equation}
As in Section~\ref{eqsym}, this implies
\begin{proposition}
If $(\Sper^\sigma)$ denotes the dual basis of $(\Mper_\sigma)$, the graded
dual $\SGSym := \SGQSym^*$ is free over the set
\begin{equation}
\{ \Sper^\alpha \,|\, \alpha \text{\ connected} \}.
\end{equation}
\end{proposition}

Indeed, Equation~(\ref{coprodM}) is equivalent to
\begin{equation}
\Sper^\alpha \Sper^\beta = \Sper^{\alpha\sconc\beta}.
\end{equation}
\qed

Note that $\SGSym$ is both a subalgebra and a quotient of $\ESym$, since
$\SGQSym$ is both a quotient and a subalgebra of $\EQSym$.

Now, as before, $\SGSym$ being a graded connected cocommutative Hopf algebra,
it follows from the Cartier-Milnor-Moore theorem that
\begin{equation}
\SGSym = U(L),
\end{equation}
where $L$ is the Lie algebra of its primitive elements.

The same argument as in Section~\ref{eqsym} proves
\begin{theorem}
The graded Lie algebra $L$ of primitive elements of $\SGSym$ is free over a
set indexed by connected permutations.
\qed
\end{theorem}

\begin{corollary}
$\SGSym$ is isomorphic to $H_O$, the Grossman-Larson Hopf algebra of
heap-ordered trees \cite{GL}.
\qed
\end{corollary}

According to~\cite{AS}, $\SGQSym$ ($=\SGSym^*$) is therefore isomorphic to the
quotient of $\FQSym$ by its coradical filtration.

%Recall that since the generating series of the number of connected
%permutations is $1+t+3t^2+13t^3+\cdots$, the Hilbert series of $L$ is
%$l(t)=t+t^2+4t^3+17t^4+\cdots$.
%
%%%%%%%%%%%%%%%%%%%%%%%%%%%%%%%%%%%%%%%%%%%%%%%%%%%%%%%%%%%%%%%%%%%%%%%%%%%%%%%
\subsection{Cyclic tensors and $\SGQSym$}

For a vector space $V$, let $\Gamma^n V$ be the subspace of $V^{\otimes n}$
spanned by \emph{cyclic tensors}, \emph{i.e.}, sums of the form
\begin{equation}
\sum_{k=0}^{n-1} (v_1\otimes \cdots\otimes v_n) \gamma^k,
\end{equation}
where $\gamma$ is the cycle $(1,2,\ldots,n)$, the right action of permutations
on tensors being as usual
\begin{equation}
(v_1\otimes \cdots\otimes v_n) \sigma =
v_{\sigma(1)}\otimes \cdots \otimes v_{\sigma(n)}.
\end{equation}

Clearly, $\Gamma^n V$ is stable under the action of $GL(V)$, and its character
is the symmetric function ``cyclic character''~\cite{ST,LST}:
\begin{equation}
l_n^{(0)} = \frac{1}{n} \sum_{d|n} \phi(d) p_d^{n/d},
\end{equation}
where $\phi$ is Euler's function.

Let $L_n^{(0)}$ be the subspace of $\CC \SG_n$ spanned by cyclic permutations.
This is a submodule of $\CC\SG_n$ for the
conjugation action
$\rho_\tau(\sigma)=\tau\sigma\tau^{-1}$ with Frobenius characteristic
$l_n^{(0)}$.
Then one can define the analytic functor $\Gamma$~\cite{Joy,Mcd}:
\begin{equation}
\Gamma (V) = \sum_{n\geq0} V^{\otimes n} \otimes_{\CC\SG_n} L_n^{(0)}.
\end{equation}

Let $\overline{\Gamma} (V) = \bigoplus_{n\geq1} \Gamma^n(V)$. Its symmetric
algebra $H(V)=S(\overline{\Gamma} (V))$ can be endowed with a Hopf algebra
structure, by declaring the elements of $\overline{\Gamma} (V)$ primitive.

As an analytic functor, $V\mapsto H(V)$ 
%S(\overline{\Gamma} (V))$
corresponds to the sequence of $\SG_n$-modules $M_n = \CC \SG_n$ endowed with
the conjugation action, that is,
\begin{equation}
%S(\overline{\Gamma} (V))
H(V)
= \bigoplus_{n\geq0} V^{\otimes n}\otimes_{\CC\SG_n} M_n,
\end{equation}
so that basis elements of $H_n(V)$ can be identified with symbols
$\qbin{w}{\sigma}$ with $w\in V^{\otimes n}$ and $\sigma\in\SG_n$ subject to
the equivalences
\begin{equation}
\qbin{w \tau^{-1}}{\tau\sigma\tau^{-1}} \equiv \qbin{w}{\sigma}.
\end{equation}
Let $A:=\{a_n\,|\, n\geq 1\}$ be an infinite linearly ordered alphabet, and
$V=\CC A$. We identify a tensor product of letters
$a_{i_1}\otimes\cdots\otimes a_{i_n}$
with the corresponding word $w=a_{i_1}\ldots a_{i_n}$ and denote by
$(w)\in\Gamma^n V$ the circular class of $w$. A basis of $H_n(V)$ is then
given by the commutative products
\begin{equation}
\underline{m} = (w_1)\cdots(w_p) \in S(\overline{\Gamma} V)
\end{equation}
of circular words, with $|w_1|+\cdots+|w_p|=n$ and $|w_i|\geq1$ for all $i$.

With such a basis element, we can associate a permutation by the following
standardization process. Fix a total order on circular words, for example the
lexicographic order on minimal representatives. Write $\underline{m}$ as a
non-decreasing product
\begin{equation}
\underline{m} = (w_1)\cdots(w_p) \text{ with } (w_1)\leq
(w_2)\leq\cdots\leq(w_p),
\end{equation}
and compute the ordinary standardization $\sigma'$ of the word
$w=w_1\cdots w_p$. Then, $\sigma$ is the permutation obtained by parenthesing
the word $\sigma'$ like $\underline{m}$ and interpreting the factors as
cycles. For example, if
\begin{equation}
\begin{split}
\underline{m} &= (cba)(aba)(ac)(ba) = (aab)(ab)(ac)(acb) \\
w &= aababacacb \\
\sigma' &= 12637495\,10\,8 \\
\sigma &= (126)(37)(49)(5\,10\,8) \\
\sigma &= (2,6,7,9,10,1,3,5,4,8) \\
\end{split}
\end{equation}
We set $\sigma = \cstd(\underline{m})$ and define it as the \emph{circular
standardized} of $\underline m$.

Let $H_\sigma(V)$ be the subspace of $H_n(V)$ spanned by those
$\underline{m}$ such that $\cstd(\underline{m})=\sigma$, and let
$\pi_\sigma : H(V)\to H_\sigma(V)$ be the projector associated with the direct
sum decomposition
\begin{equation}
H(V) = \bigoplus_{n\geq0} \bigoplus_{\sigma\in\SG_n} H_\sigma(V).
\end{equation}
Computing the convolution of such projectors then yields the following
\begin{theorem}
The $\pi_\sigma$ span a subalgebra of the convolution algebra
$\End^{\gr} H(V)$, isomorphic to $\SGQSym$ via $\pi_\sigma\mapsto M_\sigma$.
\qed
\end{theorem}

\Proof
First, note that $\pi_{\alpha}*\pi_{\beta}$ is a sum of $\pi_\gamma$:
indeed,
\begin{equation}
\pi_{\alpha}*\pi_{\beta}(\underline{m}) =
\pi_{\alpha}*\pi_{\beta}(\cstd(\underline{m})),
\end{equation}
since the circular standardized of a subword of $\underline m$ is equal to the
circular standardized of the same subword of $\cstd(\underline{m})$.
So
\begin{equation}
\pi_{\alpha}*\pi_{\beta} =  \sum_{\sigma} D_{\alpha,\beta}^\sigma \pi_\sigma.
\end{equation}
Now, by definition of $\pi_\alpha$ and by the third interpretation of the
product of the $M_\sigma$, one concludes that
$D_{\alpha,\beta}^\sigma$ is equal to the $C_{\alpha,\beta}^\sigma$ of
Equation~(\ref{defC}).
\qed

%%%%%%%%%%%%%%%%%%%%%%%%%%%%%%%%%%%%%%%%%%%%%%%%%%%%%%%%%%%%%%%%%%%%%%%%%%%%%%%
\subsection{Interpretation of $H(V)$}

The Hopf algebra $H(V)$ can be interpreted as an algebra of functions as
follows. Assume $V$ has dimension $d$, and let $a_1,\ldots,a_d$ be a basis of
$V$.

To the generator $(a_{i_1}\cdots a_{i_n})$ of $H(V)$, we associate the
function of $d$ matrices
\begin{equation}
f_{(i_1,\ldots,i_n)}(A_1,\ldots,A_d) = \tr(A_{i_1}\cdots A_{i_n}).
\end{equation}
These functions are clearly invariant under simultaneous conjugation
$A_i\mapsto MA_iM^{-1}$, and it is easy to prove that they generate the ring
of invariants of $GL(N,\CC)$ in the symmetric algebra
\begin{equation}
M_N(\CC)^{\oplus d} \simeq M_N(\CC) \otimes V.
\end{equation}
Indeed, let us set $U=\CC^N$ and identify $M_N(\CC)$ with
$U\otimes U^*$. Then, the character of $GL(U)$ in $S^n(U\otimes U^*\otimes V)$
is $h_n(XX^\vee N)$, where $X=\sum x_i$, $X^\vee=\sum x_i^{-1}$.

By the Cauchy formula,
\begin{equation}
h_n(XX^\vee N) = \sum_{\lambda\vdash n} s_\lambda(NX) s_\lambda(X^\vee),
\end{equation}
and the dimension of the invariant subspace is
\begin{equation}
\begin{split}
\dim\ S^n(U\otimes U^*\otimes V)^{GL(U)} &= \< h_n(XX^\vee N), 1\>_{GL(U)}\\
&= \sum_{\lambda\,\vdash n} \<s_\lambda(NX), s_\lambda(X) \> \\
&= \sum_{\lambda,\,\mu\,\vdash n} s_\lambda*s_\mu(N) \<s_\lambda, s_\mu\>\\
&= \sum_{\lambda\,\vdash n} (s_\lambda*s_\lambda)(N).
\end{split}
\end{equation}
The characteristic of the conjugation action of $\SG_n$ on $\CC\SG_n$ is
precisely $\sum_{\lambda\vdash n} (s_\lambda*s_\lambda)$, so this is the
dimension of $H_n(V)$.
We have therefore established:
\begin{theorem}
Let $F_N^{(d)}$ be the algebra of $GL(N,\CC)$-invariant polynomial functions
on $M_N(\CC)^{\oplus d} \simeq M_N(\CC)\otimes V$, endowed with the
comultiplication
\begin{equation}
\Delta f(A'_1,\ldots,A'_d ; A''_1,\ldots,A''_d) :=
f(A'_1\oplus A''_1,\ldots,A'_d\oplus A''_d)
\end{equation}
Then the map $(a_{i_1}\cdots a_{i_k}) \mapsto f_{(i_1,\ldots,i_k)}$ is an
epimorphism of bialgebras $H(V)\to F_N^{(d)}$.
\qed
\end{theorem}

%%%%%%%%%%%%%%%%%%%%%%%%%%%%%%%%%%%%%%%%%%%%%%%%%%%%%%%%%%%%%%%%%%%%%%%%%%%%%%%
\subsection{Subalgebras of $\SGQSym$}
\label{subsgqsym}

%%%%%%%%%%%%%%%%%%%%%%%%%%%%%%%%%%%%%%%%%%%%%%%%%%%%%%%%%%%%%%%%%%%%
\subsubsection{Symmetric functions in noncommuting variables (dual)}

For a permutation $\sigma\in\SG_n$, let $\csupp(\sigma)$ be the partition
$\pi$ of the set $[n]$ whose blocks are the supports of the cycles of
$\sigma$. The sums
\begin{equation}
\upi_{\pi} := \sum_{\csupp(\sigma)=\pi} \Mper_\sigma
\end{equation}
span a Hopf subalgebra $\PiQSym$ of $\SGQSym$, which is isomorphic to the
graded dual of the Hopf algebra of symmetric functions in noncommuting
variables (such as in~\cite{SR,BRRZ}, not to be confused with $\NCSF$), which
we will denote here by $\WSym(A)$, for \emph{Word symmetric functions}.
Indeed, from the product rule of the $\Mper_\sigma$ given in
Equation~(\ref{defC}), one easily finds
\begin{equation}
\upi_{\pi'}\upi_{\pi''} := \sum C_{\pi',\pi''}^\pi \upi_\pi,
\end{equation}
where $C_{\pi',\pi''}^\pi$ is the number of ways of splitting the parts of
$\pi$ into two subpartitions whose standardized words are $\pi'$ and $\pi''$.
%runs into the multiset of set partitions where some
%parts have as standardized word $\pi'$ and their complement $\pi''$.
For example,
\begin{equation}
\begin{split}
\upi_{ \{\{1,2,4\},\{3\}\}}\upi_{\{ \{1\}\}} &=
%\upi_{ \{\{1,2,4\},\{3\},\ \{5\} \}} +
%\upi_{ \{\{1,2,5\},\{3\},\ \{4\} \}} +
%\upi_{ \{\{1,2,5\},\{4\},\ \{3\} \}} \\
%&\ \  +
%\upi_{ \{\{1,3,5\},\{4\},\ \{2\} \}} +
%\upi_{ \{\{2,3,5\},\{4\},\ \{1\} \}} \\
%&=
\upi_{ \{\{1,2,4\},\{3\},\{5\} \}} +
2\upi_{ \{\{1,2,5\},\{3\},\{4\} \}} \\ &+
\upi_{ \{\{1,3,5\},\{4\},\{2\} \}} +
\upi_{ \{\{2,3,5\},\{4\},\{1\} \}}.
\end{split}
\end{equation}

The dual $\WSym(A)$ of $\PiQSym$ is the subspace of $\KK\<A\>$ spanned by the
orbits of $\SG(A)$ on $A^*$.
These orbits are naturally labelled by set partitions of $[n]$, the orbit
corresponding to a partition $\pi$ being constituted of the words
\begin{equation}
w = a_1\ldots a_n,
\end{equation}
such that $a_i=a_j$ iff $i$ and $j$ are in the same block of $\pi$.
The sum of these words will be denoted by $\Mw_\pi$.
%In the sequel, all set partitions will be denoted as

For example,
\begin{equation}
\Mw_{\{\{1,3,6\},\{2\},\{4,5\}\}} := \sum_{a\not=b ; b\not=c; a\not=c}
abacca.
\end{equation}

It is known that the natural coproduct of $\WSym$ (given as usual by the
ordered sum of alphabets) is cocommutative~\cite{BRRZ} and that $\WSym$ is
free over connected set partitions.
The same argument as in Theorem~\ref{prim-lib} shows that $\PiQSym^*$ is free
over the same graded set, hence that $\PiQSym$ is indeed isomorphic to
$\WSym^*$.

%%%%%%%%%%%%%%%%%%%%%%%%%%%%%%%%%%%%%%%%%
\subsubsection{Quasi-symmetric functions}

One can also embed $\QSym$ into $\PiQSym$: take as total ordering on finite
sets of integers $\{i_1<\cdots<i_r\}$ the lexicographic order on the words
$(i_1,\ldots,i_r)$. Then, any set partition $\pi$ of $[n]$ has a canonical
representative $B$ as a non-decreasing sequence of blocks $(B_1\leq B_2\leq
\cdots\leq B_r)$. Let $I=K(\pi)$ be the composition $(|B_1|,\ldots,|B_r|)$ of
$n$.
The sums
\begin{equation}
\label{def-uq}
\uq_{I} := \sum_{K(\pi)=I} \upi_\pi = \sum_{K(\sigma)=I} \Mper_\sigma
\end{equation}
where $K(\sigma)$ denotes the ordered cycle type of $\sigma$, span a Hopf
subalgebra of $\PiQSym$ and $\SGQSym$, which is isomorphic to $\QSym$.
Indeed, from the product rule of the $\Mper_\sigma$ given in
Equation~(\ref{defC}), one easily finds
\begin{equation}
\label{uq}
\uq_{I'}\uq_{I''} := \sum_{I} C_{I',I''}^I\uq_I,
\end{equation} 
where $C_{I',I''}^I$ is the coefficient of $I$ in $I'\shuffle I''$.
For example,
\begin{equation}
\begin{split}
\uq_{(1,3,1)} \uq_{(1,2)} &=
2\uq_{(1,1,2,3,1)} + 2\uq_{(1,1,3,1,2)} + 2\uq_{(1,1,3,2,1)} \\
&\ \ + \uq_{(1,2,1,3,1)} + 2\uq_{(1,3,1,1,2)} + \uq_{(1,3,1,2,1)}.
\end{split}
\end{equation}

%%%%%%%%%%%%%%%%%%%%%%%%%%%%%%%%%%%
\subsubsection{Symmetric functions}

Furthemore, if we denote by $\Lambda(I)$ the partition associated with a
composition $I$ by sorting $I$ and by $\Lambda(\pi)$ the partition $\lambda$
whose parts are the sizes of the blocks of $\pi$, the sums
\begin{equation}
\ul_{\lambda} := \sum_{\Lambda(I)=\lambda} \uq_{I}
 = \sum_{\Lambda(\pi)=\lambda} \upi_\pi
 = \sum_{Z(\sigma)=\lambda} \Mper_\sigma
\end{equation}
where $Z(\sigma)$ denotes the cycle type of $\sigma$, span a Hopf subalgebra
of $\QSym$, $\PiQSym$, and $\SGQSym$, which is isomorphic to $Sym$ (ordinary
symmetric functions).
%This last property directly comes from
%Equation~(\ref{uq}) and the embedding of $Sym$ into $\QSym$.
As an example of the product, one has
\begin{equation}
\label{ul}
\ul_{(3,3,2,1)}\ul_{(3,1,1)} = 9 \ul_{(3,3,3,2,1,1,1)}.
\end{equation}
Indeed, it follows from Equation~(\ref{uq}) that an explicit Hopf embedding of
$Sym$ into $\SGQSym$ is given by
\begin{equation}
j: p_\lambda^* \to \ul_\lambda
\end{equation}
where $p_\lambda^* = \frac{p_\lambda}{z_\lambda}$ is the adjoint basis of
products of power sums.
The images of the usual generators of $Sym$ under this embedding have simple
expressions in terms of the infinite matrix $X=(x_{ij})_{i,j\leq1}$:
\begin{equation}
j(p_n) = \tr(X^n)
\end{equation}
which implies that $j(e_n)$ is the sum of the diagonal minors of order $n$ of
$X$:
\begin{equation}
j(e_n) = \sum_{i_1<\cdots<i_n} \sum_{\sigma\in\SG_n}
         \varepsilon(\sigma) x_{i_1 i_{\sigma(1)}}\ldots x_{i_n i_{\sigma(n)}}
\end{equation}
and $j(h_n)$ is the sum of the same minors of the permanent
\begin{equation}
j(h_n) = \sum_{i_1<\cdots<i_n} \sum_{\sigma\in\SG_n}
         x_{i_1 i_{\sigma(1)}} \ldots x_{i_n i_{\sigma(n)}}.
\end{equation}
More generally, the sum of the diagonal immanants of type $\lambda$ gives
\begin{equation}
j(s_\lambda) =
\sum_{i_1<\cdots<i_n} \sum_{\sigma\in\SG_n}
\chi^\lambda(\sigma) x_{i_1 i_{\sigma(1)}} \ldots x_{i_n i_{\sigma(n)}}.
\end{equation}

\subsubsection{Involutions}

Finally, one can check that the $M_\sigma$ with $\sigma$ involutive span a
Hopf subalgebra of $\SGQSym$. Since the number of involutions of $\SG_n$ is
equal to the number of standard Young tableaux of size $n$, this algebra can
be regarded as a commutative version of the Poirier-Reutenauer
algebra~\cite{PR} denoted by $\FSym$ and realized in~\cite{NCSF6}.
Notice that this version is also isomorphic to the image of $\SGQSym$ in the
quotient of $R/{\mathcal J}$ by the relations
\begin{equation}
x_{ij}x_{jk}=0, \text{\ for all $i\not=k$.}
\end{equation}
This construction generalizes to the algebras built on permutations of
arbitrary given order.

%%%%%%%%%%%%%%%%%%%%%%%%%%%%%%%%%%%%%%%%%%%%%%%%%%%%%%%%%%%%%%%%%%%%%%%%%%%%%%%
\subsection{Quotients of $\WSym$}
\label{ws-plong}

Since we have built a subalgebra of $\PiQSym$ isomorphic to $\QSym$, we can
define an embedding
\begin{equation}
i: \QSym \hookrightarrow \PiQSym=\WSym^*,
\end{equation}
so that, dually, there is a Hopf epimorphism
$i^*: \WSym \twoheadrightarrow \NCSF$.

The dual basis $V^I$ of the $U_I$ defined in Equation~\ref{def-uq} can be
identified with equivalence classes of $\Sper^\sigma$ under the relation
\begin{equation}
\Sper^\sigma \simeq \Sper^{\sigma'}  \text{\ iff\ }
K(\sigma) = K(\sigma').
\end{equation}
The $\Sper^\sigma$ with $\sigma$ a full cycle are primitive, so that we can
take for $V_n$ any sequence of primitive generators of $\NCSF$.

It turns out that there is another natural epimorphism from $\WSym$ to
$\NCSF$.
Using the canonical ordering of set partitions introduced in
Section~\ref{subsgqsym}, that is, the lexicographic ordering on the
nondecreasing representatives of the blocks, we can as above associate a
composition $K(\pi)$ with $\pi$ and define an equivalence relation
\begin{equation}
\pi\sim\pi' \text{\ iff\ } K(\pi)=K(\pi').
\end{equation}
Then, the ideal $\goth I$ of $\WSym$ generated by the differences
\begin{equation}
\Mw_\pi - \Mw_{\pi'},\qquad \pi\sim\pi',
\end{equation}
is a Hopf ideal, and the quotient
\begin{equation}
\WSym/ {\goth I}
\end{equation}
is isomorphic to $\NCSF$. The images $V_I$ of the $\Mw_\pi$ by the canonical
projection are analogs of the monomial symmetric functions in $\NCSF$. Indeed,
the commutative image $v_\lambda$ of $\Mw_\pi$ is proportional to a monomial
function:
\begin{equation}
v_\lambda = \prod_{i} m_i(\lambda)! m_\lambda.
\end{equation}
If we introduce the coefficients $c_\lambda$ by
\begin{equation}
v_1^n = \sum_{\lambda\vdash n} c_\lambda v_\lambda
\end{equation}
then, the multivariate polynomials
\begin{equation}
B_n(x_1,\ldots,x_n) = \sum_{\lambda\vdash n} c_\lambda x_\lambda,
\end{equation}
where $x_\lambda:= x_{\lambda_1}\cdots x_{\lambda_n}$, are the exponential
Bell polynomials defined by
\begin{equation}
\sum_{n\geq0} \frac{t^n}{n!} B_n(x_1,\ldots,x_n)
= e^{\sum_{n\geq1}\frac{x_n}{n!} t^n}.
\end{equation}

%%%%%%%%%%%%%%%%%%%%%%%%%%%%%%%%%%%%%%%%%%%%%%%%%%%%%%%%%%%%%%%%%%%%%%%%%%%%%%%
\subsection{The stalactic monoid}

The constructions of Section~\ref{ws-plong} can be interpreted in terms of a
kind of Robinson-Schensted correspondence and of a plactic-like monoid.
The \emph{stalactic congruence} is the congruence $\equiv$ on $A^*$ generated
by the relations
\begin{equation}
a\,w\,a \equiv a\,a\,w,
\end{equation}
for all $a\in A$ and $w\in A^*$.

Each stalactic class has a unique representative, its \emph{canonical
representative} of the form
\begin{equation}
a_1^{m_1} a_2^{m_2} \ldots a_r^{m_r}
\end{equation}
with $a_i\not=a_j$ for $i\not=j$.

We can represent such a canonical word by a tableau-like planar diagram,
\emph{e.g.},
\begin{equation}
\PetitTableau
c^3ad^3b^2 \longleftrightarrow
\Tableau{c&a&d&b\\ c&\ &d&b \\ c&\ &d&\ \\}\,\,.
%\Tableau{c\\c\\c\\} \Tableau{a&d&b\\ c&\ &d&b \\ c&\ &d&\ \\}\,\,.
\end{equation}
Now, there is an obvious algorithm, which consists in scanning a word from
left to right and arranging its identical letters in columns, creating a new
column to the right when one scans a letter for the first time.
Let us call $P(w)$ the resulting canonical word, or, equivalently, its planar
representation.
We can compute $P(w)$ along with a $Q$-symbol recording the intermediate
shapes of the $P$-symbol. For example, to insert $w=cabccdbdd$, we have the
steps 

%\begin{equation}
%\begin{array}{ccccccc}
%\PetitTableau
%\emptyset,\emptyset
% &\fleche{c}& \Tableau{c\\}\,,\ \Tableau{1\\}
% &\fleche{a}& \Tableau{c&a\\}\,,\ \Tableau{1&2\\}
% &\fleche{b}& \Tableau{c&a&b\\}\,,\ \Tableau{1&2&3\\} 
%\end{array}
%\end{equation}

\begin{equation*}
\TasseTableau
\begin{CD}
\emptyset,\emptyset
 @>c>> \hbox{$\Tableau{c\\}\,,\ \Tableau{1\\}$}
 @>a>> \hbox{$\Tableau{c&a\\}\,,\ \Tableau{1&2\\}$}
 @>b>> \hbox{$\Tableau{c&a&b\\}\,,\ \Tableau{1&2&3\\} $}
\\[5mm]
 @>c>> \hbox{$\Tableau{c&a&b\\ c\\}\,,\ \Tableau{1&2&3\\ 4\\}$}
 @>c>> \hbox{$\Tableau{c&a&b\\ c\\ c\\}\,,\ \Tableau{1&2&3\\ 4\\ 5\\}$}
 @>d>> \hbox{$\Tableau{c&a&b&d\\ c\\ c\\}\,,\ \Tableau{1&2&3&6\\ 4\\ 5\\}$}
\\[5mm]
 @>b>> \hbox{$\Tableau{c&a&b&d\\ c&\ &b\\ c\\}\,,\ \Tableau{1&2&3&6\\ 4&\ &7\\ 5\\}$}
 @>d>> \hbox{$\Tableau{c&a&b&d\\ c&\ &b&d\\ c\\}\,,\
              \Tableau{1&2&3&6\\ 4&\ &7&8\\ 5\\}$}
 @>d>> \hbox{$\Tableau{c&a&b&d\\ c&\ &b&d\\ c&\ &\ &d\\}\,,\
              \Tableau{1&2&3&6\\ 4&\ &7&8\\ 5&\ &\ &9\\}$}
\end{CD}
\end{equation*}

Clearly, the $Q$-symbol can be interpreted as a set partition of $[n]$, whose
blocks are the columns.
In our example,
\begin{equation}
Q(cabccdbd) = \{ \{1,4,5\},\ \{2\},\ \{3,7\},\ \{6,8,9\} \}.
\end{equation}

We now see that the natural basis $\Mw_\pi$ of $\WSym$ can be characterized as
\begin{equation}
\Mw_\pi = \sum_{Q(w)=\pi} w.
\end{equation}
This is completely similar to the definition of the bases
$R_I$ of $\NCSF$ via the hypoplactic congruence,
of $\SS_t$ of $\FSym$ via the plactic congruence, and
of $\Pp_T$ of $\PBT$ via the sylvester congruence~\cite{NCSF4,NCSF6,HNT}.

Note that the defining relation of the stalactic monoid can be presented in a
plactic-like way as
\begin{equation}
a\,u\,ba \equiv a\,u\,ab
\end{equation}
for all $a,b\in A$ and $u\in A^*$.

%%%%%%%%%%%%%%%%%%%%%%%%%%%%%%%%%%%%%%%%%%%%%%%%%%%%%%%%%%%%%%%%%%%%%%%%%%%%%%%
\subsection{Other Hopf algebras derived from the stalactic congruence}

It is always interesting to investigate the behaviour of certain special
classes of words under analogs of the Robinson-Schensted correspondence.
This section presents a few examples leading to interesting combinatorial
sequences.

%%%%%%%%%%%%%%%%%%%%%%%%%%%%%%%%%
\subsubsection{Parking functions}

Recall that a word on $A=\{1,\ldots\}$ is a parking function if its
nondecreasing reordering $u_1\cdots u_k$ satisfies $u_i\leq i$.
The number of stalactic classes of parking functions of size $n$ can be
combinatorially determined as follows.

Since the congruence does not change the evaluation of a word, and since if a
word is a parking function, then all its permutations are, one can restrict
to a rearrangement class, containing a unique  nondecreasing parking
function. Those are known to be counted by Catalan numbers. Now, the
rearrangement class of a nondecreasing parking function $p$ has exactly $l!$
congruence classes if $l$ is the number of different letters of $p$.

The counting of nondecreasing parking functions $p$ by their number of
different letters is obviously the same as the counting of Dyck paths by their
number of peaks, given by the Narayana triangle (sequence
A001263 in~\cite{Slo}). To get the number of stalactic classes of parking
functions, one has to multiply the $i$th column by $i!$, that is the
definition of the unsigned Lah numbers (sequence A089231 in~\cite{Slo}),
which count with the additional parameter of number of lists, the number of
sets of lists (sequence A000262 in~\cite{Slo}), that is, the number of set
partitions of $[n]$ with an ordering inside each block but no order among the
blocks.  The first values are
\begin{equation}
1,\ 3,\ 13,\ 73,\ 501,\ 4051,\ 37633,\ 394353,\ 4596553, \ldots
\end{equation}
whereas the first rows of the Narayana and unsigned Lah triangles are
\begin{figure}[ht]
$$
\begin{array}{ccccccc}
  1 &       &      &     &     &     \\
  1 &     1 &      &     &     &     \\
  1 &     3 &    1 &     &     &     \\
  1 &     6 &    6 &   1 &     &     \\
  1 &    10 &   20 &  10 &   1 &     \\
  1 &    15 &   50 &  50 &  15 &   1 \\
\end{array}
\qquad 
\begin{array}{cccccccc}
   1 &      &     &      &     &     \\
   1 &    2 &     &      &     &     \\
   1 &    6 &   6 &      &     &     \\
   1 &   12 &  36 &   24 &     &     \\
   1 &   20 & 120 &  240 & 120 &     \\
   1 &   30 & 300 & 1200 & 1800& 720 \\
\end{array}
$$
\caption{\label{NarLah}The Narayana and unsigned Lah triangles.}
\end{figure}

One can also obtain this last result by pure algebraic calculations as
follows.
The Frobenius characteristic of the representation of $\SG_n$ on $\PF_n$ is
\begin{equation}
ch(\PF_n) = \frac{1}{n+1} h_n((n+1)X)
= \frac{1}{n+1} \sum_{\mu\vdash n} m_\mu(n+1) h_\mu(X).
\end{equation}
In this expression, each $h_\mu(X)$ is the characteristic of the permutation
representation on a rearrangement class of words, with $\mu_1$ occurrences of
some letter $i_1$, $\mu_2$ of some other letter $i_2$, and so on.
Hence, the number of stalactic classes in each such rearrangement class is
$l(\mu)!$, and the total number of stalactic classes of parking functions is
\begin{equation}
a_n = \frac{1}{n+1} \sum_{\mu\vdash n} m_\mu(n+1) l(\mu)!
\end{equation}
Since $g(z) = \sum_{n\geq0} z^n \ch(\PF_n)$ solves the functional
equation~\cite{NTlag}
\begin{equation}
\label{gg}
g(z) = \sum_{n\geq0} z^n h_n(X) g(z)^n,
\end{equation}
we see that the exponential generating function is
\begin{equation}
A(z) = \sum_{n\geq0} \alpha_n \frac{z^n}{n!} = \int_0^\infty e^{-x}g(z)dx,
\end{equation}
in the special case where $h_n=x$ for all $n\geq1$, that is
\begin{equation}
g(z) = 1 + x \frac{zg(z)}{1-zg(z)}
\end{equation}
so that
\begin{equation}
A(z) = \exp\left(\frac{z}{1-z}\right).
\end{equation}
Hence, $\alpha_n$ is the number of `sets of lists', giving back sequence
A000262 of~\cite{Slo}.
It would be interesting to find a natural bijection between the stalactic
classes of parking functions and sets of lists, compatible with the algebraic
structures.

Now, recall that the Hopf algebra of parking functions $\PQSym^*$ is spanned
by polynomials
\begin{equation}
\G_\park (A) := \sum_{\Park(w)=\park} w.
\end{equation}
As usual, it is easy to show that the equivalence defined by
\begin{equation}
G_{\park'} (A) \equiv \G_{\park''}(A)
\end{equation}
iff $\park'$ and $\park''$ are stalactically congruent, is such that the
quotient $\PQSym^*/\equiv$ is a Hopf algebra.

%%%%%%%%%%%%%%%%%%%%%%%%%%%%%
\subsubsection{Endofunctions}

The same methods allow to see that the number $\beta_n$ of stalactic classes
of endofunctions on $n$ letters is given by
\begin{equation}
\beta_n = \sum_{k=0}^n \binom{n-1}{k-1} \binom{n}{k} k!
\end{equation}
It is sequence A052852 of~\cite{Slo}, whose first terms are
\begin{equation}
1, 4, 21, 136, 1045, 9276, 93289, 1047376, 12975561, 175721140, \ldots
\end{equation}
and whose exponential generating series is
\begin{equation}
\frac{z}{1-z} e^{\frac{z}{1-z}}.
\end{equation}
As before, the classes counted by those numbers can be counted with the
additionnal parameter given by the number of different letters of their
representatives, so that one gets a new triangle, the \emph{Endt} triangle.
Moreover, this triangle is obtained by multiplying column $i$ by $i!$ from a
more classical triangle, counting the integer compositions of $2n$ in $n$
parts with a given number of ones, the \emph{Tw} triangle. Here are the
first rows of both triangles:
%f := (n,p) -> binomial(n,p) * binomial(n-1,p-1);
%g := (n,p) -> p! * f(n,p);
%h := n -> _plus(g(n,p)$p=0..n);
%
\begin{figure}[ht]
$$
\begin{array}{ccccccc}
  1 &       &       &     &      &     \\
  2 &     1 &       &     &      &     \\
  3 &     6 &     1 &     &      &     \\
  4 &    18 &    12 &   1 &      &     \\
  5 &    40 &    60 &  20 &   1  &     \\
  6 &    75 &   200 & 150 &  30  &   1 \\
\end{array}
\qquad 
\begin{array}{cccccccc}
   1 &      &      &      &      &     \\
   2 &    2 &      &      &      &     \\
   3 &   12 &    6 &      &      &     \\
   4 &   36 &   72 &   24 &      &     \\
   5 &   80 &  360 &  380 & 120  &     \\
   6 &  150 & 1200 & 3600 & 3600 & 720 \\
\end{array}
$$
\caption{\label{TwEndt}The Tw and Endt triangles.}
\end{figure}

From the algebraic point of view, the Frobenius characteristic of $[n]^n$ is
\begin{equation}
\ch([n]^n) = h_n(nX) = \sum_{\mu\vdash n} m_\mu(n) h_\mu(X),
\end{equation}
so that $\beta_n = \sum_{\mu\vdash n} m_\mu(n) l(\mu)!$, and the same method
gives directly the exponential generating series of $\beta_n$.

As above, the quotient of $\EQSym$ by $\SEF^{f'}\equiv\SEF^{f''}$ iff $f'$ and
$f''$ are stalactically congruent is a Hopf algebra.

%%%%%%%%%%%%%%%%%%%%%%%%%%%%%
\subsubsection{Initial words}

Recall that initial words are words on the alphabet of integers so that, if
letter $n$ appears in $w$, then $n-1$ also appears.
The same method allows to see that the number $\gamma_n$ of stalactic classes
of initial words on $n$ letters is
\begin{equation}
\gamma_n = \sum_{k=0}^n \binom{n-1}{k-1} k!
\end{equation}
It is sequence A001339 of~\cite{Slo}, whose first terms are
\begin{equation}
1, 3, 11, 49, 261, 1631, 11743, 95901, 876809, 8877691, 98641011, \ldots
\end{equation}
and whose exponential generating series is
\begin{equation}
\frac{e^z}{(1-z)^2}
\end{equation}
%sg := series(exp(x)/(1-x)^2,x,10);
%
As before, the classes counted by those numbers can be counted with the
additionnal parameter given by the number of different letters of their 
representatives, so that one gets a new triangle, the \emph{Arr} triangle.
Moreover, this triangle is obtained by multiplying column $i$ by $i!$ in
the Pascal triangle.
Here are the first rows of both triangles:
%f := (n,p) -> binomial(n-1,p-1);
%g := (n,p) -> p! * f(n,p);
%h := n -> _plus(g(n,p)$p=0..n);
\begin{figure}[ht]
$$
\begin{array}{ccccccc}
  1 &       &       &     &      &     \\
  1 &     1 &       &     &      &     \\
  1 &     2 &     1 &     &      &     \\
  1 &     3 &     3 &   1 &      &     \\
  1 &     4 &     6 &   4 &   1  &     \\
  1 &     5 &    10 &  10 &   5  &   1 \\
\end{array}
\qquad 
\begin{array}{cccccccc}
   1 &      &      &      &      &     \\
   1 &    2 &      &      &      &     \\
   1 &    4 &    6 &      &      &     \\
   1 &    6 &   18 &   24 &      &     \\
   1 &    8 &   36 &   96 &  120 &     \\
   1 &   10 &   60 &  240 &  600 & 720 \\
\end{array}
$$
\caption{\label{PasArr}The Pascal and Arr triangles.}
\end{figure}

Again, the stalactic quotient of $\WQSym$~\cite{NT-cras,Hiv} is a Hopf
algebra.

%%%%%%%%%%%%%%%%%%%%%%%%%%%%
\subsubsection{Generic case}

Let $A$ be an alphabet, and $V=\CC A$.
The generic symmetric function
\begin{equation}
f_n = \sum_{\mu\vdash n} l(\mu)! m_\mu(X)
\end{equation}
is the character of $GL(V)$ in the image of $V^{\otimes n}$ in the quotient of
$T(V)\sim \KK\langle A\rangle$ by the stalactic congruence.
%fn := (n) -> _plus(nops(lambda)! * sym::m(lambda)$lambda in
%combinat::partitions::list(n));
It is Schur-positive, and can be explicitly expanded on the Schur basis.
Indeed, one has
\begin{equation}
\label{fnx}
f_n = \sum_{k=0}^n c_k s_{(n-k,k)}.
\end{equation}
The coefficients $c_k$ are given by sequence A000255 of~\cite{Slo}, whose
first terms are
\begin{equation}
1, 1, 3, 11, 53, 309, 2119, 16687, 148329, 1468457, 16019531, \dots
\end{equation}
%sg := series(exp(-x)/(1-x)^2,x,10);
To see this, write again
\begin{equation}
l(\mu)! = \int_0^\infty t^{l(\mu)} e^{-t}dt.
\end{equation}
Then
\begin{equation}
f = \sum_{n\geq0} f_n
= \int_0^\infty e^{-t} \sum_\mu t^{l(\mu)}m_\mu dt.
\end{equation}
and
\begin{equation}
\label{devel}
\begin{split}
\sum_{\mu} t^{l(\mu)}m_\mu = \prod_{i\geq1} \left(1+\frac{tx_i}{1-x_i}\right)
&= \prod_{i\geq1} \frac{1-(1-t)x_i}{1-x_i}\\
%&= \sigma_1((1-q)X)|_{q=1-t}\\
%&= \lambda_{t-1}(X)\sigma_1(X)\\
&= \left(\sum_{k\geq0}(t-1)^ke_k(X) \right)\left( \sum_{l\geq0}h_l(X)\right).
\end{split}
\end{equation}

Since
\begin{equation}
\int_0^\infty e^{-t} (t-1)^kdt=  d_k,
\end{equation}
the number of {derangements} in $\SG_n$, we have finally
\begin{equation}
f = \sum_{n\geq0} \sum_{k=0}^n d_k e_k(X) h_{n-k}(X).
\end{equation}
Expanding
\begin{equation}
e_k h_{n-k} = s_{(n-k,1^k)} + s_{n-k+1,1^{k-1}},
\end{equation}
we get Formula~(\ref{fnx}).
Alternatively, we can express $c_k$ as 
\begin{equation}
c_k = \int_{0}^\infty e^{-t} t (t-1)^k dt
\end{equation}
%int( e^(-t)*t*(t-1)^k,t,0..infinity)
since the term of degree $n$ in Equation~(\ref{devel}) is
\begin{equation}
t \sum_{k=0}^n (t-1)^k s_{(n-k,1^k)}.
\end{equation}
The exponential generating series of these numbers is given by
\begin{equation}
\frac{e^{-z}}{(1-z)^2}.
\end{equation}

%%%%%%%%%%%%%%%%%%%%%%%%%%%%%%%%%%%%%%%%%%%%%%%%%%%%%%%%%%%%%%%%%%%%%%%%%%%%%%%
%%%%%%%%%%%%%%%%%%%%%%%%%%%%%%%%%%%%%%%%%%%%%%%%%%%%%%%%%%%%%%%%%%%%%%%%%%%%%%%
%%%%%%%%%%%%%%%%%%%%%%%%%%%%%%%%%%%%%%%%%%%%%%%%%%%%%%%%%%%%%%%%%%%%%%%%%%%%%%%
\section{Structure of $\SGSym$}

%%%%%%%%%%%%%%%%%%%%%%%%%%%%%%%%%%%%%%%%%%%%%%%%%%%%%%%%%%%%%%%%%%%%%%%%%%%%%%%
\subsection{A realization of $\SGSym$}

In the previous section, we have built a commutative algebra of permutations
$\SGQSym$ from explicit polynomials on a set of auxiliary variables
$x_{i\,j}$.
One may ask whether its non-commutative dual admits a similar realization in
terms of non-commuting variables $a_{i\,j}$.

We shall find such a realization, in a somewhat indirect way, by first
building from scratch a Hopf algebra of permutations
$\PhiSym \subset \KK\,\langle\, a_{i\,j} \,|\, i,j\geq 1\,\rangle$,
whose operations can be described in terms of the cycle structure of
permutations.
Its coproduct turns out to be cocommutative, and the isomorphism with
$\SGSym$ follows as above from the Cartier-Milnor-Moore theorem.

Let $\{a_{i\,j}, i,j\geq1\}$ be an infinite set of non-commuting
indeterminates. We use the biword notation
\begin{equation}
a_{i\,j} \equiv \ncbinomial{i}{j}, \quad
\ncbinomial{i_1}{j_1} \cdots \ncbinomial{i_n}{j_n} \equiv
\ncbinomial{i_1\ldots i_n}{j_1\ldots j_n}
\end{equation}

Let $\sigma\in\SG_n$ and let $(c_1,\ldots,c_k)$ be a decomposition of $\sigma$
into disjoint cycles. With any cycle, one associates its \emph{cycle words},
that is, the words obtained by reading the successive images of any
element of the cycle.
For example, the cycle words associated with the cycle $(1342)$ are $1342,
2134, 3421, 4213$.
For all word $w$, we denote by $C(w)$ the \emph{cycle} associated with the
inverse of its standardized word.

We now define
\begin{equation}
\myphi_\sigma := \sum_{x,a} \ncbinomial{x}{a},
\end{equation}
where the sum runs over all words $x$ such that $x_i=x_j$ iff $i$ and
$j$ belong to the same cycle of $\sigma$, and such that the standardized word
of the subword $x$ of $a$ whose letter positions belong to cycle $c_l$
satisfied $C(x)=c_l$.
By extension, such biwords will be said to have $\sigma$ as \emph{cycle
decomposition}.
%the inverse of the standardized word of a cycle word of $c_l$.

Note that any biword appears in the expansion of exactly one $\myphi_\sigma$.

\begin{example}
\begin{equation}
\myphi_{12} = \sum_{x\not= y} \ncbinomial{x\ y}{a\ b}.
\end{equation}
\begin{equation}
\myphi_{41352} =
  \sum_{x\not=y ; \Std(abde)^{-1}=1342,\, 2134,\, 3421, \text{\rm\ or\ } 4213}
  {\ncbinomial{x\ x\ y\ x\ x}{a\ b\ c\ d\ e}}.
\end{equation}
\end{example}

\begin{theorem}
The $\myphi_\sigma$ span a subalgebra $\PhiSym$ of
$\KK\,\langle\, a_{ij} \,|\, i,j\geq1 \,\rangle$.
More precisely,
%there exist non-negative integers
%$g_{\alpha,\beta}^{\sigma}$ $(0$ or $1)$ such that
\begin{equation}
\myphi_\alpha \myphi_\beta = \sum g_{\alpha,\beta}^{\sigma} \myphi_\sigma,
\end{equation}
where $g_{\alpha,\beta}^{\sigma}\in \{0,1\}$.
%\medskip
%
Moreover, $\PhiSym$ is free over the set
\begin{equation}
\{\myphi_\alpha \,|\, \alpha \text{\ connected} \}.
\end{equation}
\end{theorem}

\Proof
Let $\sigma'$ and $\sigma''$ be two permutations and let $w$ be a biword
appearing in $\myphi_{\sigma'} \myphi_{\sigma''}$.
The multiplicity of $w$ is 1. To get the first part of the theorem, we only
need to prove that all words $w'$ appearing in the same $\myphi_\sigma$ as $w$
also appear in this product. Given a letter $x$ appearing in the first
row of $w$ at some positions, the subwords $a$ of the elements of
$\myphi_\sigma$ taken from the \emph{second} row at those positions, have the
same image by $C$ as the corresponding element of $w$.

Thus, we onyl have to prove that all words $w$ having the same
image by $C$ satisfy that all their prefixes (and suffixes) of a given length
have also same image by $C$. It is sufficient to prove the result on
permutations. Now, given two permutations $\sigma$ and $\tau$, $\sigma^{-1}$
and $\tau^{-1}$ are cycle words of the same cycle, iff for some $k$,
$\tau=\gamma^k\sigma$ where $\gamma$ is the cyclic permutation
$(12\cdots n)$.
Since is also the case for the standardized words of the prefixes of both
permutations of a given length, the property holds.

Moreover, any biword can be uniquely written as a concatenation of a maximal
number of biwords such that no letter appears in the first row of two
different biwords and that the letters of the second row of a biword are all
smaller than the letters of the second row of the next one. This proves
that the $\myphi_\alpha$ where $\alpha$ is connected are free. The usual
generating series argument then proves that those elements generate $\PhiSym$.
\qed

To give the precise expression of the product $\myphi_\alpha\myphi_\beta$, we 
first need to define two operations on cycles.

The first operation is just the circular shuffle on disjoint cycles: if $c'_1$
and $c''_1$ are two disjoint cycles, their \emph{cyclic shuffle}
$c'_1 \picyc c''_1$ is the set of cycles $c_1$ such that their cycle words are
obtained by applying the usual shuffle on the cycle words of $c'_1$ and
$c''_1$. This definition makes sense because a shuffle of cycle words
associated with two words on disjoint alphabets splits as a union of cyclic
classes.

For example, the cyclic shuffle $(132)\picyc(45)$ gives the set of cycles
\begin{equation}
\begin{split}
\{ (13245), (13425), (13452), (14325), (14352), (14532), \\
   (13254), (13524), (13542), (15324), (15342), (15432) \}.
\end{split}
\end{equation}
These cycles correspond to the following list of permutations which are
those appearing in Equation~(\ref{312-54b}), except for the first one which
will be found later:
\begin{equation}
\begin{split}
\{ &\ 34251,\ 35421,\ 31452,\ 45231,\ 41532,\ 41253,  \\
   &\ 35214,\ 34512,\ 31524,\ 54213,\ 51423,\ 51234 \}.
\end{split}
\end{equation}

Let us now define an operation on two sets $C_1$ and $C_2$ of disjoint cycles.
We call \emph{matching} a list of all those cycles, some of the cycles being
paired, always one of $C_1$ with one of $C_2$. The cycles remaining alone are
considered to be associated with the empty cycle. With all matchings associate
the set of sets of cycles obtained by the cyclic product $\picyc$ of any pair
of cycles.
The union of those sets of cycles is denoted by $C_1\picycgen C_2$.

For example, the matchings corresponding to
$C_1=\{(1),(2)\}$ and $C_2=\{(3),(4)\}$ are:
\begin{equation}
\label{c1c2a}
\begin{split}
&\{(1)\}\{(2)\}\{(3)\}\{(4)\},\
\{(1)\}\{(2),(3)\}\{(4)\},\ 
\{(1)\}\{(2),(4)\}\{(3)\}, \\
&\{(1),(3)\}\{(2)\}\{(4)\},\
\{(1),(3)\}\{(2),(4)\},\
\{(1),(4)\}\{(2)\}\{(3)\},\
\{(1),(4)\}\{(2),(3)\},
\end{split}
\end{equation}
and the product $C_1\picycgen C_2$ is then
\begin{equation}
\label{c1c2b}
\begin{split}
&\{(1),(2),(3),(4)\},\ \
\{(1),(23),(4)\},\ \
\{(1),(24),(3)\}, \\
&\{(13),(2),(4)\},\ \
\{(13),(24)\},\ \
\{(14),(2),(3)\},\ \
\{(14),(23)\}.
\end{split}
\end{equation}

Remark that this calculation is identical with the Wick formula in quantum
field theory (see~\cite{BO} for an explanation of this coincidence).

We are now in a position to describe the product $\myphi_\sigma \myphi_\tau$:

\begin{proposition}
\label{cycdec}
let $C_1$ be the cycle decomposition of $\sigma$ and $C_2$ be the cycle
decomposition of $\tau$, shifted by the size of $\sigma$. Then the
permutations indexing the elements appearing in the product
$\myphi_\sigma \myphi_\tau$ are the permutations whose cycle decompositions
belong to $C_1\picycgen C_2$.
\end{proposition}

\Proof
Recall that the product $\myphi_\sigma \myphi_\tau$ is a sum of
biwords with multiplicity one since a word appears in exactly one
$\myphi_\sigma$. So we only have to prove that biwords appearing in
$\myphi_\sigma \myphi_\tau$ are the same as biwords whose cycle decomposition
are in $C_1\picycgen C_2$.
First, by definition of the cyclic shuffle and of operation~$\picycgen$, if a
biword has its cycle decomposition in $C_1\picycgen C_2$, its
prefix of size $n$ has cycle decomposition $C_1$ whereas its suffix of
size $p$ has cycle decomposition $C_2$, where $n$ (resp. $p$) is the size of
$C_1$ (resp. $C_2$).

Conversely, let $w_1$ (resp. $w_2$) be a biword with cycle decomposition $C_1$
(resp. $C_2$) and let us consider $w=w_1\cdot w_2$. For all letters in the
first row of $w$, either it only appears in $w_1$, either only in $w_2$, or in
both $w_1$ and $w_2$. In the first two cases, we obtain the corresponding
cycle of $C_1$ (or $C_2$, shifted). In the last case, the cycle decomposition
of the word of the second row corresponding to this letter belongs to the
cyclic shuffle of the corresponding cycles of $C_1$ and $C_2$ (hence matching
those two cycles). Indeed, if $\Std^{-1}(w_1)$ is a cycle word of $c_1$ and
$\Std^{-1}(w_2)$ is a cycle word of $c_2$, then $\Std^{-1}(w)$ is a cycle
word of an element of the cyclic shuffle $c_1\shuffle c_2$: it is already
known when $\Std^{-1}(w_1)=c_1$ and $\Std^{-1}(w_2)=c_2$ since it is the
definition of $\FQSym$~\cite{NCSF6} and cyclying $c_1$ or $c_2$ only amounts
to cycle $c_2$, hence copying the definition of the cyclic shuffle.
\qed

For example, with $\sigma=\tau=12$, one finds that $C_1=\{(1),(2)\}$ and
$C_2=\{(3),(4)\}$. It is then easy to check that one goes from
Equation~(\ref{c1c2b}) to Equation~(\ref{12-12b}) by computing the
corresponding permutations.

\begin{example}
{\rm
\begin{equation}
\label{12-43b}
\myphi_{12} \myphi_{21} = \myphi_{1243} + \myphi_{1342} + \myphi_{1423} +
\myphi_{3241} + \myphi_{4213}.
\end{equation}
\begin{equation}
\label{12-12b}
\myphi_{12} \myphi_{12} = \myphi_{1234} + \myphi_{1324} + \myphi_{1432} +
\myphi_{3214} + \myphi_{3412} + \myphi_{4231} + \myphi_{4321}.
\end{equation}
\begin{equation}
\myphi_{1} \myphi_{4312} = \myphi_{15423} + \myphi_{25413} + \myphi_{35421} +
\myphi_{45123} + \myphi_{51423}.
\end{equation}
\begin{equation}
\begin{split}
\label{312-54b}
\myphi_{312} \myphi_{21} &=\ \ \
\myphi_{31254} + \myphi_{31452} + \myphi_{31524} + \myphi_{34251} +
\myphi_{34512} + \myphi_{35214} + \myphi_{35421} \\
&\ \ +\ \myphi_{41253} + \myphi_{41532} +  \myphi_{45231} + \myphi_{51234} +
\myphi_{51423} + \myphi_{54213}.
\end{split}
\end{equation}
}
\end{example}

Let us recall a rather general recipe to obtain the coproduct of a
combinatorial Hopf algebra from a realization in terms of words on an ordered
alphabet $X$.
Assume that $X$ is the ordered sum of two mutually commuting alphabets
$X'$ and $X''$.  Then define the coproduct as $\Delta(F)=F(X'\dot{+}X'')$,
identifying $F'\tensor F''$ with $F'(X')F''(X'')$~\cite{NCSF6, NTpfb}.

There are many different ways to define a coproduct on $\PhiSym$ compatible
with the realization since there are many ways to order an alphabet of
biletters: order the letters of the first alphabet, order the letters of the
second alphabet, or order lexicographically with respect to one alphabet and
then to the second.

In the sequel, we only consider the coproduct obtained by ordering the
biletters with respect to the first alphabet so that $\myphi_\sigma$ is
primitive if $\sigma$ consists of only one cycle.
More precisely, thanks to the definition of the $\myphi$, it is easy to see
that it corresponds to the unshuffling of the cycles of a permutation:

\begin{equation}
\Delta\myphi_\sigma := \sum_{(\alpha,\beta)}
\myphi_\alpha \tensor \myphi_\beta,
\end{equation}
where the sum is taken over all pairs of permutations $(\alpha,\beta)$ such
that $\alpha$ is obtained by renumbering the elements of a subset of cycles
of $\sigma$ (preserving the relative order of values), and $\beta$ by doing
the same on the complementary subset of cycles.
For example, if $\sigma=(1592)(36)(4)(78)$, the subset $(1592)(4)$
gives $\alpha=(1452)(3)$ and $\beta=(12)(34)$.

\begin{example}
\begin{equation}
\Delta\myphi_{12} = \myphi_{12}\tensor 1 + 2 \myphi_{1}\tensor\myphi_{1}
+1 \tensor \myphi_{12}.
\end{equation}
\begin{equation}
\Delta\myphi_{312} = \myphi_{312}\tensor 1 + 1 \tensor \myphi_{312}.
\end{equation}
\begin{equation}
\Delta\myphi_{4231} = \myphi_{4231}\tensor 1 + 2 \myphi_{321}\tensor\myphi_{1}
+ \myphi_{21}\tensor \myphi_{12} + \myphi_{12}\tensor \myphi_{21}
+ 2 \myphi_{1}\tensor \myphi_{321} + 1\tensor \myphi_{4231}.
\end{equation}
\end{example}

The following theorem is a direct consequence of the definition of the
coproduct on the realization.
% It is also a consequence of a more general result presented in~\cite{NTbar}.

\begin{theorem}
$\Delta$ is an algebra morphism, so that $\PhiSym$ is a graded
bialgebra (for the grading $\deg\myphi_\sigma = n$ if $\sigma\in\SG_n$).
Moreover, $\Delta$ is cocommutative.
\qed
\end{theorem}

The same reasoning as in Section~\ref{sgqsym} shows that

\begin{theorem}
$\SGSym$ and $\PhiSym$ are isomorphic as Hopf algebras.
\qed
\end{theorem}

To get the explicit isomorphism sending $\myphi$ to $\Sper$, let us
first recall that a \emph{connected permutation} is a permutation $\sigma$
such that $\sigma([1,k])\not = [1,k]$ for any $k\in [1,n-1]$. Any permutation
$\sigma$ has a unique maximal factorization
$\sigma=\sigma_1\bullet\cdots\bullet\sigma_r$ into connected permutations.
We then define
\begin{equation}
\label{sprime}
\Sper'_\sigma := \myphi_{\sigma_1}\cdots\myphi_{\sigma_r}.
\end{equation}

Then
\begin{equation}
\Sper'_\sigma = \myphi_{\sigma_1\bullet\cdots\bullet\sigma_r} +
\sum_{\mu}\myphi_\mu
\end{equation}
where all permutations $\mu$ have strictly less cycles than $r$.
So the $\Sper'$ form a basis of $\PhiSym$.
Moreover, they are a multiplicative basis with product given by shifted
concatenation of permutations, so that they multiply as the $\Sper$ do.
Moreover, the coproduct of $\Sper'_\sigma$ is the same as for
$\myphi_\sigma$, so the same as for $\Sper^\sigma$.
So both bases $\Sper$ and $\Sper'$ have the same product and the same
coproduct.
This proves

\begin{proposition}
The linear map $ \Sper^\sigma \mapsto \Sper'_\sigma $
realizes the Hopf isomorphism between $\SGSym$ and $\PhiSym$.
\qed
\end{proposition}

There is another natural isomorphism: define
\begin{equation}
\label{ssec}
\Sper''_\sigma := \sum_{x,a} \ncbinomial{x}{a},
\end{equation}
where the sum runs over all words $x$ such that $x_i=x_j$ if (but \emph{not}
only if) $i$ and $j$ belong the same cycle of $\sigma$ and such that the
standardized word of the subword of $a$ consisting of the indices of cycle
$c_l$ is equal to the inverse of the standardized word of a cycle word of
$c_l$.

That both bases $\Sper'$ and $\Sper''$ have the same product and
coproduct comes from the fact that if $(c_1)\cdots(c_p)$ is the cycle
decomposition of $\sigma$,

\begin{equation}
\Sper''_\sigma = \sum_{(c);(c)\in
                 (c_1)\picycgen(c_2)\picycgen\cdots\picycgen(c_p)}
                     \myphi_{(c)}.
\end{equation}
For example,
\begin{equation}
\begin{split}
\Sper''_{2431} = \Sper''_{(124)(3)}  & =
\myphi_{(124)(3)} + \myphi_{(1423)} + \myphi_{(1234)} + \myphi_{(1324)}\\
& =\myphi_{2431}  + \myphi_{4312}   + \myphi_{2341}   + \myphi_{3421}.
\end{split}
\end{equation}

%%%%%%%%%%%%%%%%%%%%%%%%%%%%%%%%%%%%%%%%%%%%%%%%%%%%%%%%%%%%%%%%%%%%%%%%%%%%%%%
\subsection{Quotients of $\PhiSym$}

Let $I$ be the ideal of $\PhiSym$ generated by the differences
\begin{equation}
\phi_{\sigma} - \phi_{\tau}
\end{equation}
where $\sigma$ and $\tau$ have the same cycle type.

The definitions of its product and coproduct directly imply that $I$ is a Hopf
ideal.
Since the cycle types are parametrized by integer partitions, the quotient
$\PhiSym/I$ has a basis $Y_\lambda$, corresponding to the class of
$\phi_{\sigma}$, where $\sigma$ has $\lambda$ as cycle type.

From Equations~(\ref{12-43b})-(\ref{312-54b}), one finds:

\begin{example}
\begin{equation}
Y_{11} Y_{2} = Y_{211} + 4 Y_{31}, \qquad
Y_{11}^2 = Y_{1111} + 2 Y_{22} + 4 Y_{211}.
\end{equation}
\begin{equation}
Y_1 Y_4 = Y_{41} + 4 Y_5, \qquad Y_{3} Y_{2} = Y_{32} + 12 Y_5.
\end{equation}
\end{example}

\begin{theorem}
$\PhiSym/I$ is isomorphic to $Sym$, the Hopf algebra of ordinary symmetric
functions,
\medskip

If one writes
$\lambda=(\lambda_1,\ldots,\lambda_p)=(1^{m_1},\ldots,k^{m_k})$,
an explicit isomorphism is given by
\begin{equation}
Y_\lambda \mapsto 
\frac{\prod_{i=1}^{k} m_i!}{\prod_{j=1}^p (\lambda_j-1)!} m_\lambda.
\end{equation}
\end{theorem}

\Proof
This follows from the description of
$\myphi_\sigma \myphi_\tau$ given in Proposition~\ref{cycdec}.
\qed

%%%%%%%%%%%%%%%%%%%%%%%%%%%%%%%%%%%%%%%%%%%%%%%%%%%%%%%%%%%%%%%%%%%%%%%%%%%%%%%
%%%%%%%%%%%%%%%%%%%%%%%%%%%%%%%%%%%%%%%%%%%%%%%%%%%%%%%%%%%%%%%%%%%%%%%%%%%%%%%
%%%%%%%%%%%%%%%%%%%%%%%%%%%%%%%%%%%%%%%%%%%%%%%%%%%%%%%%%%%%%%%%%%%%%%%%%%%%%%%
\section{Parking functions and trees}

%%%%%%%%%%%%%%%%%%%%%%%%%%%%%%%%%%%%%%%%%%%%%%%%%%%%%%%%%%%%%%%%%%%%%%%%%%%%%%%
\subsection{A commutative algebra of parking functions}

It is also possible to build a commutative pendant of the Hopf algebra of
parking functions introduced in~\cite{NT1}:
let $\PF_n$ be the set of parking functions of length $n$.
For $\park\in\PF_n$, set, as before

\begin{equation}
\Mpa_\park := \sum_{i_1<\cdots < i_n} x_{i_1\ i_{\park(1)}} \cdots
x_{i_n\ i_{\park(n)}}.
\end{equation}

Then, using once more the same arguments of Section 2, we conclude that the
$\Mpa_\park$ form a linear basis of a $\ZZ$-subalgebra $\CPQSym$ of $\EQSym$,
which is also a sub-coalgebra if one defines the coproduct in the usual way.
%that is, from special cuts in graphs (see~\cite{NTT} for more details).

\begin{example}
\begin{equation}
\Mpa_{1} \Mpa_{11} = \Mpa_{122} + \Mpa_{121} + \Mpa_{113}.
\end{equation}
\begin{equation}
\Mpa_{1} \Mpa_{221} = \Mpa_{1332} + \Mpa_{3231} + \Mpa_{2231} + \Mpa_{2214}.
\end{equation}
\begin{equation}
\Mpa_{12} \Mpa_{21} = \Mpa_{1243} + \Mpa_{1432} + \Mpa_{4231} + \Mpa_{1324} +
\Mpa_{3214} + \Mpa_{2134}.
\end{equation}
\begin{equation}
\Delta \Mpa_{525124}  = \Mpa_{525124} \tensor 1 + 1\tensor \Mpa_{525124}.
\end{equation}
\begin{equation}
\Delta \Mpa_{4131166} = \Mpa_{4131166} \tensor 1 +
                        \Mpa_{41311}\tensor \Mpa_{11} +
                        1 \tensor \Mpa_{4131166}.
\end{equation}
\end{example}

The main interest of the non-commutative and non-cocommutative Hopf algebra of
parking functions defined in~\cite{NT1} is that it leads to two algebras of
trees. The authors obtain a cocommutative Hopf algebra of planar binary trees
by summing over the distinct permutations of parking functions, and an algebra
of planar trees by summing over hypoplactic classes.

We shall now investigate whether similar constructions can be found for the
commutative version $\CPQSym$.

%%%%%%%%%%%%%%%%%%%%%%%%%%%%%%%%%%%%%%%%%%%%%%%%%%%%%%%%%%%%%%%%%%%%%%%%%%%%%%%
\subsection{From labelled to unlabelled parking graphs}

A first construction, which can always be done for Hopf algebras of labelled
graphs is to build a subalgebra by summing over labellings.
Notice that this subalgebra is the same as the subalgebra obtained by summing
endofunctions graphs over their labellings. Those graphs are also known as
endofunctions (hence considered there as unlabelled graphs) in~\cite{BLL}.

The dimension of this subalgebra in degree $n$ is equal to the number of
unlabelled parking graphs
\begin{equation}
1, 1, 3, 7, 19, 47, 130, 343, 951, 2615, 7318, 20491, 57903, 163898,\ldots
\end{equation}
known as sequence A001372 in~\cite{Slo}.
For example, here are the $7$ unlabelled parking graphs of size $3$ (to be
compared with the $16$ parking functions):

\begin{figure}[ht]
\epsfig{file=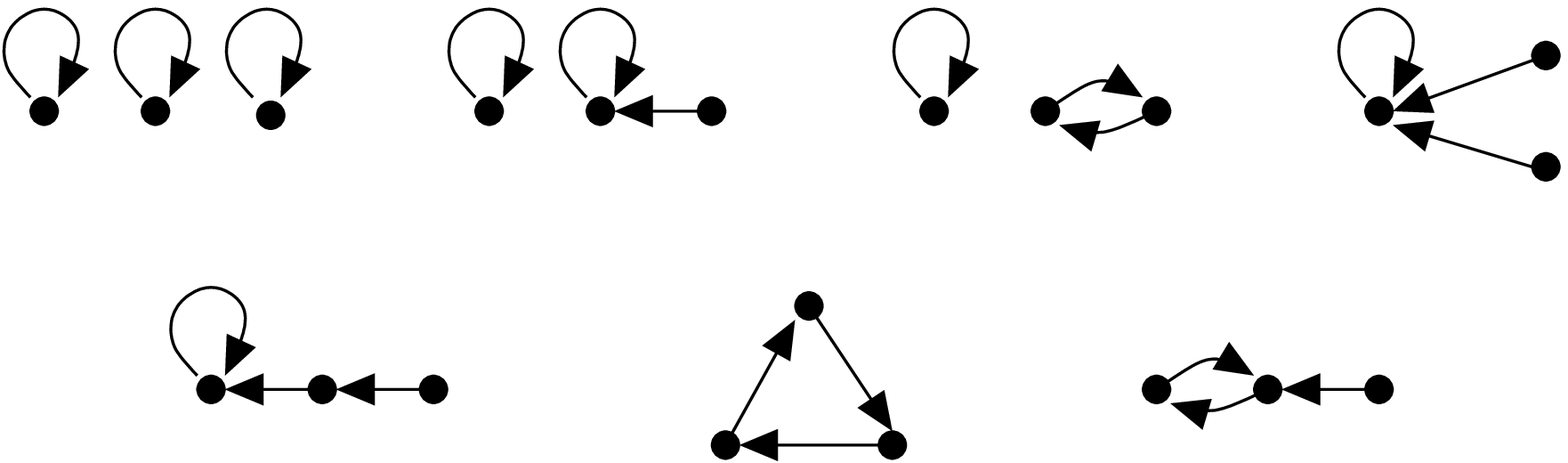,width=12cm}
\end{figure}

The product of two such unlabelled graphs is the concatenation of graphs and
the coproduct of an unlabelled graph is the unshuffle of its connected
subgraphs. So this algebra is isomorphic to the polynomial algebra on
generators indexed by connected parking graphs.

%%%%%%%%%%%%%%%%%%%%%%%%%%%%%%%%%%%%%%%%%%%%%%%%%%%%%%%%%%%%%%%%%%%%%%%%%%%%%%%
\subsection{Binary trees and nondecreasing parking functions}

One can easily check that in $\CPQSym$, summing over parking functions having
the same reordering does not lead to a subalgebra. However, if we denote by
$I$ the subspace of $\CPQSym$ spanned by the $\Mpa_\park$ where $\park$ is not
nondecreasing, $I$ is an ideal and a coideal, and $\CCQSym := \CPQSym/I$ is
therefore a commutative Hopf algebra with basis given by the classes
$\Mpa_{\parkc}:= \overline{\Mpa_{\park}}$
labelled by nondecreasing parking functions.

Notice that $\CCQSym$ is also isomorphic to the image of $\CPQSym$ in the
quotient of $R/{\mathcal J}$ by the relations
\begin{equation}
x_{ij}x_{kl}=0, \text{\ for all $i<k$ and $j>l$.}
\end{equation}

The dual basis of $\Mpa_{\parkc}$ is
\begin{equation}
S^{\parkc} := \sum_{\park} S^\park,
\end{equation}
where the sum is taken over all permutations of $\parkc$.

The dual $\CCQSym^*$ is free over the set $S^\parkc$, where $\parkc$ runs over
connected nondecreasing parking functions. So if one denotes by $\CQSym$ the
Catalan algebra defined in~\cite{NT1}, 
the usual Cartier-Milnor-Moore argument then shows that
\begin{equation}
\CQSym \sim \CCQSym^* ,\qquad\qquad \CCQSym \sim \CQSym^*.
\end{equation}

%%%%%%%%%%%%%%%%%%%%%%%%%%%%%%%%%%%%%%%%%%%%%%%%%%%%%%%%%%%%%%%%%%%%%%%%%%%%%%%
\subsection{From nondecreasing parking functions to rooted forests}

Nondecreasing parking functions correspond to parking graphs of a particular
type: namely, rooted forests with a particular labelling (it corresponds to
nondecreasing maps), the root being given by the loops in each connected
component.

Taking sums over the allowed labellings of a given rooted forest, we get that
the
\begin{equation}
M_F := \sum_{supp(\parkc)=F} M_\parkc,
\end{equation}
span a commutative Hopf algebra of rooted forests, which is likely to coincide
with the quotient of the Connes-Kreimer algebra~\cite{CK} by its coradical
filtration~\cite{AS}.

%%%%%%%%%%%%%%%%%%%%%%%%%%%%%%%%%%%%%%%%%%%%%%%%%%%%%%%%%%%%%%%%%%%%%%%%%%%%%%%
%%%%%%%%%%%%%%%%%%%%%%%%%%%%%%%%%%%%%%%%%%%%%%%%%%%%%%%%%%%%%%%%%%%%%%%%%%%%%%%
%%%%%%%%%%%%%%%%%%%%%%%%%%%%%%%%%%%%%%%%%%%%%%%%%%%%%%%%%%%%%%%%%%%%%%%%%%%%%%%
\section{Quantum versions}

%%%%%%%%%%%%%%%%%%%%%%%%%%%%%%%%%%%%%%%%%%%%%%%%%%%%%%%%%%%%%%%%%%%%%%%%%%%%%%%
\subsection{Quantum quasi-symmetric functions}
\label{qsymq}

When several Hopf algebra structures can be defined on the same class of
combinatorial objects, it is tempting to try to interpolate between them.
This can be done, for example with compositions: the algebra of quantum
quasi-symmetric functions $\QSym_q$~\cite{TU,NCSF6} interpolates between
quasi-symmetric functions and non-commutative symmetric functions.

However, the natural structure on $\QSym_q$ is not exactly that of a Hopf
algebra but rather of a \emph{twisted Hopf algebra}~\cite{LZ}.

Recall that the coproduct of $\QSym(X)$ amounts to replace $X$ by the ordered
sum $X' \dot{+} X''$ of two isomorphic and mutually commuting alphabets. On
the other hand, $\QSym_q$ can be realized by means of an alphabet of
$q$-commuting letters
\begin{equation}
x_jx_i = q x_i x_j, \text{ for } j>i.
\end{equation}
Hence, if we define a coproduct on $\QSym_q$ by
\begin{equation}
\Delta_q f(X) = f(X'\dot+X''),
\end{equation}
with $X'$ and $X''$ $q$-commuting with each other, it will be an algebra
morphism
\begin{equation}
\QSym_q \to \QSym_q(X'\dot+X'') \simeq \QSym_q \otimes_{\chi} \QSym_q
\end{equation}
for the \emph{twisted product of tensors}
\begin{equation}
(a\otimes b) \cdot (a'\otimes b') = \chi(b,a') (aa' \otimes bb'),
\end{equation}
where
\begin{equation}
\label{chi}
\chi(b,a') = q^{\deg(b)\cdot \deg(a')}
\end{equation}
for homogeneous elements $b$ and $a'$.

It is easily checked that $\Delta_q$ is actually given by the same formula as
the usual coproduct of $\QSym$, that is
\begin{equation}
\Delta_q M_I = \sum_{H\cdot K=I} M_H\otimes M_K.
\end{equation}

The dual twisted Hopf algebra, denoted by $\NCSF_q$, is isomorphic to $\NCSF$
as a algebra. If we denote by $S^I$ the dual basis of $M_I$, $S^I
S^J=S^{I\cdot J}$, and $\NCSF_q$ is freely generated by the $S^{(n)}=S_n$,
whose coproduct is
\begin{equation}
\Delta_q S_n = \sum_{i+j=n} q^{ij} S_i\otimes S_j.
\end{equation}

As above, $\Delta_q$ is an algebra morphism
\begin{equation}
\NCSF_q \to \NCSF_q \otimes_{\chi} \NCSF_q,
\end{equation}
where $\chi$ is again defined by Equation~(\ref{chi}).

%%%%%%%%%%%%%%%%%%%%%%%%%%%%%%%%%%%%%%%%%%%%%%%%%%%%%%%%%%%%%%%%%%%%%%%%%%%%%%%
\subsection{Quantum free quasi-symmetric functions}

The previous constructions can be lifted to $\FQSym$.
Recall from~\cite{NCSF6} that $\phi(\F_\sigma) = q^{l(\sigma)} F_{c(\sigma)}$
is an algebra homomorphism $\FQSym \to \QSym_q$, which is in fact induced by
the specialization $\phi(a_i)=x_i$ of the underlying free variables $a_i$ to
$q$-commuting variables $x_i$.

The coproduct of $\FQSym$ is also defined by
\begin{equation}
\Delta \F(A) = \F(A' \dot+ A''),
\end{equation}
where $A'$ and $A''$ are two mutually commuting copies of $A$~\cite{NCSF6}.
If instead one sets $a''a' = qa' a''$, one obtains again a twisted Hopf
algebra structure $\FQSym_q$ on $\FQSym$, for which $\phi$ is a
homomorphism.
% The following results directly come from the description of
%$QSym_q$ recalled in the previous paragraph.
With these definitions at hand, one can see that the arguments given
in~\cite{TU} to establish the results recalled in Section~\ref{qsymq} prove in
fact the following more general result:

\begin{theorem}
~
Let $A'$ and $A''$ be $q$-commuting copies of the ordered alphabet $A$,
\emph{i.e.}, $a''a'=qa'a''$ for $a'\in A'$ and $a''\in A''$.
Then, the coproduct
\begin{equation}
\Delta_q f = f(A'\dot+ A'')
\end{equation}
defines a twisted Hopf algebra structure.
It is explictly given in the basis $\F_\sigma$ by
\begin{equation}
\label{deltaq}
\Delta_q \F_\sigma = \sum_{\alpha\cdot\beta=\sigma}
                     q^{\inv(\alpha,\beta)}
                     \F_{\Std(\alpha)}\otimes \F_{\Std(\beta)}
\end{equation}
where $\inv(\alpha,\beta)$ is the number of inversions of $\sigma$ with one
element in $\alpha$ and the other in $\beta$.
\medskip

More precisely, $\Delta_q$ is an algebra morphism with values in the
twisted tensor product of graded algebras $\FQSym\otimes_\chi \FQSym$
where $(a\otimes b)(a'\otimes b') = \chi(b,a')
(aa'\otimes bb')$ and $\chi(b,a')=q^{\deg(b).\deg(a')}$ for homogeneous
elements $b$, $a'$.
\medskip

The map $\phi: \FQSym_q \to \QSym_q$ defined by
\begin{equation}
\phi(\F_\sigma) = q^{l(\sigma)}  F_{c(\sigma)}
\end{equation}
is a morphism of twisted Hopf algebras.
\qed
\end{theorem}

\begin{example}
\begin{equation}
\Delta_q \F_{2431}  = \F_{2431}\otimes 1 + q^3 \F_{132}\otimes \F_{1}
+ q^3 \F_{12}\otimes \F_{21} + q\F_1\otimes\F_{321} + 1\otimes \F_{2431}.
\end{equation}
\begin{equation}
\Delta_q \F_{3421} =  \F_{3421}\otimes 1 + q^3 \F_{231}\otimes \F_{1}
+ q^4 \F_{12}\otimes \F_{21} + q^2\F_1\otimes\F_{321} + 1\otimes \F_{3421}.
\end{equation}
\begin{equation}
\Delta_q \F_{21} = \F_{21}\otimes 1 + q \F_{1}\otimes\F_{1} + 1\otimes\F_{21}.
\end{equation}
\begin{equation}
\begin{split}
(\Delta_q \F_{21})(\Delta_q\F_{1}) =
(\F_{213}+\F_{231}+\F_{321}) \otimes 1
+ ( \F_{21} + q^2(\F_{12}+\F_{21}))\otimes \F_{1} \\
+ \F_{1} \otimes ( q^2\F_{21} + q(\F_{12}+\F_{21}))
+1\otimes (\F_{213}+\F_{231}+\F_{321}).
\end{split}
\end{equation}
\begin{equation}
\Delta_q \F_{213} = \F_{213}\otimes 1 + \F_{21}\otimes \F_{1}
+ q \F_{1}\otimes \F_{12} + 1\otimes \F_{213}.
\end{equation}
\begin{equation}
\Delta_q \F_{231} = \F_{231}\otimes 1 + q^2 \F_{12}\otimes \F_{1}
+ q \F_{1}\otimes \F_{21} + 1\otimes \F_{231}.
\end{equation}
\begin{equation}
\Delta_q \F_{321} = \F_{321}\otimes 1 + q^2 \F_{21}\otimes \F_{1}
+ q^2 \F_{1}\otimes \F_{21} + 1\otimes \F_{321}.
\end{equation}
\end{example}

Finally, one can also define a one-parameter family of ordinary Hopf algebra
structures on $\FQSym$, by restricting formula~(\ref{deltaq}) for $\Delta_q$
to connected permutations $\sigma$, and requiring that $\Delta_q$ be an
algebra homomorphism. Then, for $q=0$, $\Delta_q$ becomes cocommutative, and
it is easily shown that the resulting Hopf algebra is isomorphic to $\SGSym$.

However, it follows from~\cite{NCSF6} that for generic $q$, the Hopf algebras
defined in this way are all isomorphic to $\FQSym$.
This suggests to interpret $\FQSym$ as a kind of quantum group: it would be
the generic element of a quantum deformation of the enveloping algebra
$\SGSym=U(L)$.
Similar considerations apply to various examples, in particular to the
Loday-Ronco algebra $\PBT$, whose commutative version obtained in~\cite{NT1}
can be quantized in the same way as $\QSym$, by means of $q$-commuting
variables~\cite{NTpfb}.

There is another way to obtain $\QSym$ from $\FQSym$: it is
known~\cite{NCSF6} that $\QSym$ is isomorphic to the image of
$\FQSym(A)$ in the hypoplactic algebra $\KK[A^*/\equiv_H]$.
One may then ask whether there exists a $q$-analogue of the hypoplactic
congruence leading directly to $\QSym_q$.

Recall that the hypoplactic congruence can be presented as the bi-sylvester
congruence:
\begin{equation}
\begin{split}
bvca \equiv bvac, \text{ with } a<b\leq c, \\
cavb \equiv acvb, \text{ with } a\leq b<c,
\end{split}
\end{equation}
and $v\in A^*$.

A natural $q$-analogue, compatible with the above $q$-commutation is
\begin{equation}
\begin{split}
bvca \equiv_{qH}\ \ q\ bvac, \text{ with } a<b\leq c, \\
cavb \equiv_{qH}\ \ q\ acvb, \text{ with } a\leq b<c,
\end{split}
\end{equation}
and $v\in A^*$.
Then, we have

\begin{theorem}
The image of $\FQSym(A)$ under the natural projection
$\KK\langle A\rangle \to \KK\langle A\rangle/\equiv_{qH}$ is isomorphic to
$\QSym_q$ as an algebra, and also as a twisted Hopf algebra for the coproduct
$A\to A'\dot+ A''$, $A'$ and $A''$ being $q$-commuting alphabets.
\end{theorem}

Moreover, it is known that if one only considers the sylvester congruence
\begin{equation}
cavb \equiv_S acvb,
\end{equation}
the quotient $\FQSym(A)$ under the natural projection
$\KK\langle A\rangle \to \KK\langle A\rangle/\equiv_S$ is isomorphic to the
Hopf algebra of planar binary trees of Loday and Ronco~\cite{LR1,HNT}. The
previous construction provides natural twisted $q$-analogs of this Hopf
algebra. Indeed, let the $q$-sylvester congruence $\equiv_{qS}$ be
\begin{equation}
cavb \equiv_{qS}\ \ q\ acvb, \text{ with } a\leq b<c.
\end{equation}
Then, since this congruence is compatible with the $q$-commutation,
we have

\begin{theorem}
The image of $\FQSym(A)$ under the natural projection
$\KK\langle A\rangle \to \KK\langle A\rangle/\equiv_{qS}$ is a twisted Hopf
algebra, with basis indexed by planar binary trees.
\qed
\end{theorem}

%%%%%%%%%%%%%%%%%%%%%%%%%%%%%%%%%%%%%%%%%%%%%%%%%%%%%%%%%%%%%%%%%%%%%%%%%%%%%%%
%%%%%%%%%%%%%%%%%%%%%%%%%%%%%%%%%%%%%%%%%%%%%%%%%%%%%%%%%%%%%%%%%%%%%%%%%%%%%%%
%%%%%%%%%%%%%%%%%%%%%%%%%%%%%%%%%%%%%%%%%%%%%%%%%%%%%%%%%%%%%%%%%%%%%%%%%%%%%%%
\section*{Acknowledgements}
This project has been partially supported by EC's IHRP Programme, grant
HPRN-CT-2001-00272, ``Algebraic Combinatorics in Europe''.
The authors would also like to thank the contributors of the MuPAD project,
and especially of the combinat part, for providing the development environment
for this research.

%%%%%%%%%%%%%%%%%%%%%%%%%%%%%%%%%%%%%%%%%%%%%%%%%%%%%%%%%%%%%%%%%%%%%%%%%%%%%%%
%%%%%%%%%%%%%%%%%%%%%%%%%%%%%%%%%%%%%%%%%%%%%%%%%%%%%%%%%%%%%%%%%%%%%%%%%%%%%%%
%%%%%%%%%%%%%%%%%%%%%%%%%%%%%%%%%%%%%%%%%%%%%%%%%%%%%%%%%%%%%%%%%%%%%%%%%%%%%%%

\end{document}